\documentclass[a4paper,12pt,nosumlimits]{amsart}
\usepackage{latexsym,amsmath,amssymb,color,mathrsfs,enumerate}
\usepackage[english]{babel}

\numberwithin{equation}{section}
\usepackage{color}

\newtheorem{thm}{Theorem}[section]
\newtheorem{cor}[thm]{Corollary}
\newtheorem{lem}[thm]{Lemma}
\newtheorem{prop}[thm]{Proposition}

\theoremstyle{definition}
\newtheorem{defn}[thm]{Definition}

\newcommand{\nn}{ \nonumber }

\newcommand{\sx}{\mathsf{x}}
\newcommand{\st}{\mathsf{z}}

\newcommand{\E}{\mathcal{E}}
\newcommand{\A}{\mathcal{A}}

\newcommand{\la}{\langle}
\newcommand{\ra}{\rangle}
\newcommand{\tr}{\triangleright}
\newcommand{\tl}{\triangleleft}

\newcommand{\ep}{\epsilon}
\newcommand{\M}{\textup{M}}
\newcommand{\SL}{\textup{SL}}
\newcommand{\U}{\textup{U}}
\newcommand{\SU}{\textup{SU}}

\newcommand{\ZZ}{\mathbb{Z}}
\newcommand{\C}{\mathbb{C}}
\newcommand{\RR}{\mathbb{R}}
\newcommand{\R}{\mathbb{R}}
\newcommand{\iq}{q^{-1}}
\newcommand{\usu}{\mathcal{U}_q(\mathfrak{su}(2))}
\newcommand{\ASU}{\A[\SU_q(2)]}
\newcommand{\p}{\partial}

\newcommand{\D}{\textup{d}}
\newcommand{\vu}{|j,m;\uparrow\ra}
\newcommand{\vd}{|j,m;\downarrow\ra}

\renewcommand{\o}{{}_{\scriptscriptstyle(1)}}
\renewcommand{\t}{{}_{\scriptscriptstyle(2)}}
\newcommand{\thr}{{}_{\scriptscriptstyle(3)}}

\title[Spin Geometry of the Quantum Two-Sphere]
{The 3D Spin Geometry of the Quantum Two-Sphere}
\author{Simon Brain}\address{\flushleft Dipartimento di Matematica e
Informatica, Universit\`{a} di Trieste, Via A. Valerio 12/1, 34127
Trieste, Italia} \email{brain@sissa.it}
\author{Giovanni Landi}\address{\flushleft Dipartimento di Matematica e
Informatica, Universit\`{a} di Trieste, Via A. Valerio 12/1, 34127
Trieste, Italia, and INFN, Sezione di Trieste, Trieste, Italia}
\email{landi@univ.trieste.it}

\date{v1: 10 March 2010; v2: 30 September 10}

\begin{document}

\begin{abstract} We study a three-dimensional differential calculus $\Omega^1S^2_q$ on the
standard Podle\'{s} quantum two-sphere $S^2_q$, coming from the
Woronowicz 4D${}_+$ differential calculus on the quantum group
$\SU_q(2)$. We use a frame bundle approach to give an explicit
description of $\Omega^1S^2_q$ and its associated spin geometry in
terms of a natural spectral triple over $S^2_q$. We equip this
spectral triple with a real structure for which the commutant
property and the first order condition are satisfied up to
infinitesimals of arbitrary order.\end{abstract}

\maketitle

\tableofcontents

\section{Introduction}

The standard quantum two-sphere $S^2_q$ has proven to be one of the
most important and useful examples in trying to understand the
relationship between the geometric/analytic world of noncommutative geometry
and the algebraic setting of quantum group theory. At the algebraic
level, it is known that $S^2_q$ has a unique left-covariant
two-dimensional differential calculus \cite{pod2,pod3}. On the other
hand, it is known that this same calculus is recovered {\em via}
analytic techniques by means of a noncommutative spin geometry
\cite{ds:est,sw04}. This compatibility has led to the discovery of other
noncommutative two-dimensional geometries on $S^2_q$ with a range of
interesting properties \cite{dalw}. In this paper, we extend the
investigation to the noncommutative spin geometry of a differential
calculus on $S^2_q$ whose dimension is equal to three.

Quantum two-spheres were constructed and classified by Podle\'{s} in
\cite{pod1}. The standard sphere $S^2_q$ is unique amongst the
Podle\'{s} family in that it also appears as the base space of the
noncommutative Hopf fibration $\SU_q(2)\to S^2_q$ constructed in
\cite{brmaj} as a basic example of a quantum principal bundle. By
equipping the total space $\SU_q(2)$ with the 3D differential
calculus of \cite{wor}, one finds that the two-dimensional
differential calculus on $S^2_q$ appears as an associated vector
bundle. This `quantum frame bundle' approach to noncommutative
geometry, developed in \cite{Ma1,Ma:sphere}, has been applied
successfully to study a host of examples, not least the
two-dimensional geometry of the quantum sphere $S^2_q$ itself.

The present paper also uses the frame bundle approach to study the
geometry of $S^2_q$, but this time starting with the
4D${}_+$ differential calculus on $\SU_q(2)$ of \cite{wor}. This
calculus has the advantage of being bicovariant under both left and
right translation, in contrast with the 3D calculus, which is only
left-covariant. Using the framing theory we recover the
three-dimensional differential calculus $\Omega^1S^2_q$ of \cite{pod2,durd1,durd2} on
$S^2_q$. The methods we use are
well-adapted to the principal bundle structure and as a consequence
we immediately find an explicit description of the bimodule
relations in $\Omega^1S^2_q$, including a decomposition into
irreducible components. We do not discuss the deeper aspects of the
Riemannian geometry such as Hodge structure and connection theory:
these will be developed elsewhere \cite{lz}.

Our main results concern the spin geometry of the three-dimensional
calculus $\Omega^1S^2_q$. Remarkably, we find that the spinor bundle
of $S^2_q$ is unchanged from the one used in
\cite{ds:est,Ma:sphere,sw04} for the two-dimensional calculus. We
construct a Dirac operator $D$ which implements the exterior
derivative in $\Omega^1S^2_q$, finding that the eigenvalues of $|D|$
grow not faster than $q^{-2j}$ for large $j$ and hence that
the associated spectral triple has metric dimension zero.

Moreover, we equip this spectral triple with a $\ZZ_2$-grading
operator and a real structure which is defined `up to compact
operators', in the sense that the `commutant property' and the
`first order condition' for a real spectral triple \cite{ac:fncg}
are satisfied up to infinitesimals of arbitrary order.
As we shall see, this is in contrast with \cite{ds:est},
where a `true' real structure for the `two-dimensional' calculus on
$S^2_q$ was given ({\em cf}. also \cite{sw04}), but is parallel to
the results of \cite{dalw} for the sphere $S^2_q$. We also find
that the `KO-theoretic' dimension of this real spectral triple is
equal to the classical value, just two.

The paper is organised as follows. In \S\ref{se:prelims} we give a
brief overview of the construction of quantum differential calculi
on quantum groups and their homogeneous spaces, followed by the
general quantum frame bundle construction itself. Following this,
\S\ref{se:pod} recalls the elementary geometry of the Hopf fibration
$\SU_q(2)\to S^2_q$ and the Hopf algebra $\usu$ which describes its
symmetries. In \S\ref{se:diffs} we describe the differential
structure of the Hopf fibration. We start from the 4D quantum
differential calculus on the total space $\SU_q(2)$ from which we
derive the calculus on the bundle fibre $\U(1)$. The structure of
the calculus $\Omega^1S^2_q$ is then obtained as a `framed quantum
manifold' in the sense of \cite{Ma:sphere}. Finally, in
\S\ref{se:spectral} we construct our spectral triple
$(\A[S^2_q],\mathcal{H},D)$ over $S^2_q$, which in addition we equip
with a $\mathbb{Z}_2$-grading $\Gamma$ of the spinor bundle
$\mathcal{H}$ and a real structure $J:\mathcal{H}\to\mathcal{H}$.

\subsubsection*{Notation}
In this paper we make frequent use of the `$q$-numbers' defined by
\begin{equation}\label{qno}[x]:=\frac{q^x-q^{-x}}{q-\iq}\end{equation}
for each $x\in \R$ and $q\not= 1$. Furthermore,
for the sake of brevity we introduce the constants
\begin{equation}\label{const}
\mu:=q+\iq,\qquad \nu:=q-\iq\end{equation} to be used throughout the
paper. Our convention is that $\mathbb{N}=\{0,1,2,\ldots\}$.

\section{Preliminaries on Quantum Principal Bundles}\label{se:prelims}

We start with some generalities on differential calculi and quantum
principal bundles. These will be endowed both with universal and
non-universal compatible calculi.

\subsection{Differential structures}\label{diffs}
Let $P$ be a complex $*$-algebra with unit. A {\em first order differential calculus} over
$P$ is a pair $(\Omega^1P,\D)$ where $\Omega^1P$ is a
$P$-$P$-bimodule (the one-forms) and $\D:P
\rightarrow \Omega^1P$ is a linear map obeying the Leibniz rule
$$\D(ab)=a(\D b)+(\D a)b,\qquad a,b\in P,$$ and such that the map $P\otimes P \rightarrow \Omega^1P$
defined by $a
\otimes b \mapsto a\D b$ is surjective.

The {\em universal differential calculus} over $P$ is the pair
$(\widetilde \Omega^1 P, \tilde \D)$, where $\widetilde \Omega^1
P:=\ker m$ is the kernel of the product map $m: P\otimes P \to P$ on $P$,
with obvious bimodule structure
$$
p\cdot(a\otimes b)=pa\otimes b,\quad (a\otimes b)\cdot p=a\otimes
bp,\qquad a,b,p\in P$$
and $\tilde \D$ is defined by $\tilde \D p:=1\otimes p-p\otimes 1$, for each $p \in P$. It is so-called
because any other differential calculus $(\Omega^1P, \D)$ over $P$
arises as a quotient $\Omega^1 P = \widetilde \Omega^1 P/N_P$, where
$N_P$ is some $P$-$P$-sub-bimodule of $\widetilde \Omega^1P$. With
the projection $\pi_{P}: \widetilde \Omega^1P \rightarrow \Omega^1P$
one has $\D = \pi_{P} \circ \tilde \D$.

If $H$ is a Hopf algebra, we write $m_H:H\otimes H\to H$ and $1_H$
for its product and unit, $\Delta_H:H\to H\otimes H$ and $\ep_H:H\to
\C$ for its coproduct and counit and $S_H:H\to H$ for its antipode
(when there is no possibility of confusion, we omit the subscript
$H$). We use Sweedler notation $\Delta(h)=h\o\otimes h\t$ for the
coproduct. A differential calculus $\Omega^1H$ over a Hopf algebra
$H$ is said to be {\em left-covariant} if the coproduct $\Delta$,
viewed as a left coaction of $H$ on itself, extends to a left
coaction $\Delta_L:\Omega^1H \rightarrow H\otimes \Omega^1H$ such
that $\D$ is an intertwiner and $\Delta_L$ is a bimodule map:
$$\Delta_L(\D h)=(\textup{id} \otimes \D) \Delta_L(h), \qquad\Delta_L(h\omega)=\Delta(h)\cdot\Delta_L(\omega), \qquad \Delta_L(\omega h)=\Delta_L(\omega)\cdot\Delta(h)$$ for all $h \in H$, $\omega \in \Omega^1
H$. A similar definition holds for a right-covariant calculus, now
with a right coaction $\Delta_R:\Omega^1H\rightarrow
\Omega^1H\otimes H$. A calculus is said to be {\em bicovariant} if
it is both left and right covariant with commuting coactions.
The universal calculus $\widetilde\Omega^1H$ is bicovariant
when equipped with the left and right tensor product coactions on $H\otimes H$.

Left-covariant differential calculi on a Hopf algebra $H$ are
classified as follows after \cite{wor}. First, it may be shown that
the linear map
\begin{equation}\label{map r}r:H \otimes H \rightarrow H \otimes H,
\qquad r(a \otimes b):=ab\o\otimes b\t,\end{equation} is an
isomorphism with inverse
\begin{equation}\label{r inv}r^{-1}:H \otimes H \rightarrow H
\otimes H,\qquad r^{-1}(a \otimes b)=aS(b\o) \otimes
b\t.\end{equation} Upon restricting $r$ to the universal calculus
$\widetilde \Omega^1 H$ we obtain an isomorphism
$$r: \widetilde \Omega^1 H \rightarrow H \otimes H^+,$$ where $H^+:=\ker
\ep_H$ denotes the augmentation ideal of $H$. This is in fact an
isomorphism of $H$-$H$ bimodules if we equip $H\otimes H^+$ with the
bimodule structure \begin{equation}\label{right module calculus}a
\cdot (b \otimes \omega)=ab \otimes \omega, \qquad (a \otimes
\omega) \cdot b = ab\o \otimes \omega b\t,\qquad a,b\in H,~\omega\in
H^+\end{equation} and an isomorphism of $H$-$H$-bicomodules if we
equip $H\otimes H^+$ with the bicomodule structure $$\Delta_L(a
\otimes \omega)=a\o \otimes \big(a\t \otimes \omega\big),\qquad
\Delta_R(a\otimes\omega)=\big(a\o\otimes\omega\o\big)\otimes a\t
\omega\t,\qquad a\in H,~\omega\in H^+.$$ Any left-covariant
sub-bimodule $N_H$ of $\widetilde \Omega^1H$ is carried to a right
ideal $I_H$ of $H^+$ by the map $r$ in~\eqref{map r}.
Conversely, any right ideal $I_H$ arises in this way from a
left-covariant sub-bimodule of $\widetilde\Omega^1H$. It follows
that the left-covariant differential calculi on $H$ are in
one-to-one correspondence with right ideals $I_H \subset H^+$;
indeed, given such an $I_H$, one has $\Omega^1H\simeq H\otimes
\Lambda^1$, where $\Lambda^1\cong H^+/I_H$ are the left-invariant
one-forms. We also write
$\Omega^1_{\textup{inv}}H:=r^{-1}(\Lambda^1)$.

A left-covariant sub-bimodule $N_H$ is also right-covariant if and
only if the corresponding ideal $I_H$ is stable under the right
adjoint coaction $$\textup{Ad}_R:H\to H\otimes H,\qquad
\textup{Ad}_R(a)=a\t\otimes S(a\o)a\thr,$$ in the sense that
$\textup{Ad}_R(I_H)\subset I_H\otimes H$. It follows that
bicovariant calculi on $H$ are in one-to-one correspondence with
right ideals $I_H$ of $H^+$ which are $\textup{Ad}_R$-stable
\cite{wor}.

Given a left-covariant differential calculus $\Omega^1H$  over $H$,
the {\em quantum tangent space} of $\Omega^1H$ is the vector space
\begin{equation}\mathcal{T}_H:=\left\{ X\in H'~|~X(1)=0 ~\text{and}~
X(a)=0~\text{for all}~a\in I_H\right\},\end{equation} where the
vector space $H'$ is the linear dual of $H$. This tangent space
admits many properties analogous to the classical case, in
particular there exists a unique bilinear form $\la
\,\cdot\,|\,\cdot\,\ra: \mathcal{T}_H\times \Omega^1H\to \C$ such
that
\begin{equation}
\la X|a\D b\ra=\ep_H(a)X(b),\qquad a,b\in
H,~X\in\mathcal{T}_H.
\end{equation}
With respect to this bilinear
form, the vector spaces $\Omega^1_{\textup{inv}}H$ and
$\mathcal{T}_H$ are non-degenerately paired, so that
$$\textrm{dim}\,\Omega^1_{\textup{inv}}H=\textrm{dim}\,\mathcal{T}_H=\textrm{dim}\,\Lambda^1.$$
This number is said to be the {\em dimension} of the left-covariant
differential calculus $\Omega^1H$.

\subsection{Quantum principal bundles}\label{qpb}The general set-up for a principal fibration of
noncommutative spaces is an algebra $P$ (playing the role of the
algebra of functions on the total space) which is a right comodule
algebra for a Hopf algebra $H$ with coaction $\delta_R:P \rightarrow
P \otimes H$. The algebra of functions on the base space of the
fibration is the subalgebra $M$ of $P$ consisting of coinvariant
elements under $\delta_R$,
$$M:=P^H=\{p \in P: \delta_R(p)=p \otimes 1\}.$$ For a well-defined
bundle structure at the level of universal differential calculi, one
requires exactness of the following sequence \cite{brmaj},
\begin{equation}\label{eqn universal exact seq}0 \rightarrow
P(\widetilde \Omega^1M)P \xrightarrow{j} \widetilde \Omega^1P
\xrightarrow{\textup{ver}} P \otimes H^+ \rightarrow
0,\end{equation} with $H^+$ the augmentation ideal, as before. The algebra
inclusion $M \hookrightarrow P$ extends to an inclusion $\tilde
\Omega^1M \hookrightarrow \widetilde \Omega^1P$ of universal
differential calculi, hence $P(\widetilde \Omega^1M)P$ are the analogues of
the horizontal one-forms (classically this corresponds to the space
of one-forms which have been pulled back from the base of the
fibration). The map $\textup{ver}$ is defined by $$\textup{ver}(p
\otimes p')=p\delta_R(p');$$ the generator of the vertical one-forms.
We say that the inclusion $M \hookrightarrow P$ is a {\em
quantum principal bundle with universal calculi} and structure
quantum group $H$. Requiring exactness of the sequence \eqref{eqn
universal exact seq} is equivalent to requiring that the induced
canonical map
\begin{equation}\label{can map}\chi:P\otimes_MP\to P\otimes H,
\qquad p\otimes_M p'\mapsto p\delta_R(p')\end{equation} be
bijective. If this is the case, one also says that the triple
$(P,H,M)$ is an $H$-{\em Hopf-Galois extension}. This bijection
condition is enough for a principal bundle structure at the level of
universal differential calculi.

For a principal bundle with non-universal calculi extra conditions are required that we briefly recall.
Assume then that $P$ and $M$ are equipped with differential calculi $\Omega^1P=\widetilde \Omega^1P/N_P$ and $\Omega^1M=\widetilde\Omega^1M/N_M$, where $N_P$ and $M_M$ are
sub-bimodules of $\widetilde \Omega^1P$ and $\widetilde \Omega^1M$
respectively. Assume further that $H$ is equipped with a
left-covariant calculus $\Omega^1H$ corresponding to a right ideal
$I_H$.

Compatibility of the differential structures means that the
calculi satisfy the conditions
\begin{equation}
N_M = N_P \cap \widetilde \Omega^1M \qquad \textup{and} \qquad
\delta_R(N_P)\subset N_P \otimes H .  \label{submodules}
\end{equation}
The role of the first condition is to ensure that $\Omega^1M$ is spanned by
elements of the form $m \D n$ with $m,n \in M$ and is hence obtained
by restricting the calculus on $P$. The second condition in \eqref{submodules}
is sufficient to ensure covariance of
$\Omega^1P$. Finally, we need the sequence
\begin{equation}\label{eqn exact seq non-univ}0 \rightarrow
P(\Omega^1M)P \rightarrow \Omega^1P \xrightarrow{\textup{ver}}P
\otimes \Lambda^1 \rightarrow 0\end{equation} to be exact. This
sequence is the analogue of the sequence (\ref{eqn universal exact
seq}) but now at the level of non-universal calculi. The
$P$-$P$-bimodule $P(\Omega^1M)P$ once again makes up the
horizontal one-forms and $\textup{ver}(p \otimes p')=p\delta_R(p')$
is the canonical map which generates the vertical one-forms. The
condition
\begin{equation}
\label{projected submod}
\textup{ver}(N_P)=P \otimes I_H
\end{equation} ensures that the map
$$\textup{ver}:\Omega^1P \rightarrow P \otimes \Lambda^1, \quad
\Lambda^1 \simeq H^+/I_H$$ is well-defined and yields that the
sequence \eqref{eqn exact seq non-univ} is indeed exact.

\subsection{Framed quantum manifolds}\label{se:frames}Suppose that the total space $P$ of the bundle is itself a
Hopf algebra equipped with a Hopf algebra surjection $\pi:P
\rightarrow H$. Here we have a coaction of $H$ on $P$ by coproduct
and projection to $H$,
$$\delta_R:P \rightarrow P \otimes H,\qquad \delta_R=(\textup{id}
\otimes \pi)\Delta.$$ The base is then the quantum homogeneous space
$M=P^H$ of coinvariants and the algebra inclusion $M\hookrightarrow
P$ is automatically an $H$-Hopf-Galois extension, {\em i.e.} a
quantum principal bundle with universal calculi. To impose
non-universal differential structure we suppose that $\Omega^1P$ is
left-covariant for $P$ and $\Omega^1H$ is left-covariant for $H$, so
that they are defined by right ideals $I_P$ and $I_H$ of $P^+$ and
$H^+$ respectively. We ensure the first of \eqref{submodules} by
taking it as a definition of $\Omega^1M$; in the case at hand, the
remaining compatibility conditions in
\eqref{submodules}--\eqref{projected submod} reduce to
\begin{equation} \label{hom bundle} (\textup{id} \otimes
\pi)\textup{Ad}_R(I_P) \subset I_P \otimes H, \qquad
\pi(I_P)=I_H.\end{equation} Thus a choice of left-covariant calculus
on $P$ satisfying these conditions automatically gives a principal
bundle with non-universal calculi \cite{Ma:sphere}.

We say that an algebra $M$ is a {\em framed quantum manifold} if it
is the base of a quantum principal bundle, $M=P^H$, to which
$\Omega^1M$ is an associated vector bundle. To give $M$ as a framed
quantum manifold we therefore require not only a quantum principal
bundle $\delta_R:P \rightarrow P \otimes H$ as above but also a
right $H$-comodule $V$, so that $\E:=(P \otimes V)^H$ plays the role
of the sections of the corresponding associated vector bundle (the space $P \otimes V$
is equipped with the tensor product coaction). Moreover, we require
a `soldering form' $\theta:V \rightarrow P\Omega^1M$ such that the
map
$$s_\theta: \E \rightarrow \Omega^1M, \quad p\otimes v \mapsto
p\theta (v)$$ is an isomorphism.

For a general $M$ it is usually not obvious how to go about looking
for a framing. However in the case of a quantum homogeneous space
with compatible calculi one has a `standard' framing in the
following way \cite{Ma:sphere}. If the conditions in~(\ref{hom bundle}) are
satisfied then the algebra $M=P^H$ is automatically framed by the
bundle $(P,H,M)$. The $H$-comodule $V$ and soldering form $\theta$
are given explicitly by the formul{\ae}
\begin{equation}\label{frame}V=(P^+ \cap M)/(I_P \cap M), \quad
\Delta_Rv= \tilde v\t \otimes S\pi(\tilde v\o),
\quad\theta(v)=S\tilde v\o \D \tilde v\t,\end{equation} with
$\tilde v$ any representative of $v$ in $P^+ \cap M$ and
$\Delta(\tilde v)=\tilde v\o\otimes \tilde v\t$ is the coproduct on $P$.

\section{The Standard Podle\'{s} Sphere}\label{se:pod}
We recall here some of the basic geometry of the so-called standard
Podle\'{s} quantum two-sphere $S^2_q$ of \cite{pod1}. We begin with
the quantum group $\ASU$ and its symmetries $\usu$, from which we
obtain the quantum sphere $S^2_q$ as the base space of the quantum
Hopf fibration $\SU_q(2)\to S^2_q$. Finally we sketch the
construction of a family of quantum line bundles over $S^2_q$ which
shall prove useful in what is to follow.

\subsection{The quantum group $\SU_q(2)$} Recall that the coordinate algebra $\A[\M_q(2)]$ of functions on the quantum matrices $\M_q(2)$ is the associative unital algebra
generated by the entries of the matrix
$$\mathbf{x}=(\sx_i{}^j)=\begin{pmatrix}a&b\\c&d\end{pmatrix}$$ obeying the relations
\begin{eqnarray}\label{mat rels}& ab=qba,\quad
ac=qca,\quad bd=qdb,\quad cd=qdc,\\\nonumber & bc=cb,\quad
ad-da=(q-\iq)bc,\end{eqnarray} with $0\ne q\in \C$ a deformation
parameter. The algebra $\A[\M_q(2)]$ has a coalgebra structure given
by $\Delta(\sx_i{}^j)=\sx_i{}^\mu\otimes\sx_\mu{}^j$ and
$\ep(x_i{}^j)=\delta_i{}^j$. From $\A[\M_q(2)]$ we obtain a Hopf
algebra $\A[\SL_q(2)]$ upon quotienting by the determinant relation
$ad=1+qbc$ (equivalently $da=1+\iq bc$) and defining an antipode by
$$S\begin{pmatrix}a&b\\c&d\end{pmatrix}=\begin{pmatrix}d&-\iq b\\-q c&a\end{pmatrix}.$$

When the deformation parameter $q$ is taken to be real $\A[\M_q(2)]$ is made into a $*$-algebra by
defining the anti-linear involution
\begin{equation}\label{star}\mathbf{x}^*=\begin{pmatrix}a^*&b^*\\c^*&d^*\end{pmatrix}:=\begin{pmatrix}d&-q c\\-\iq
b&a\end{pmatrix}.\end{equation} It is not difficult to see that
$\A[\SL_q(2)]$ inherits this $*$-structure. Without loss of
generality we take $0<q<1$. The compact quantum group $\A[\SU_q(2)]$
is defined to be the quotient of $\A[\SL_q(2)]$ by the additional
relations $S(\mathbf{x}_k{}^l)=(\mathbf{x}_l{}^k)^*$.
Thus in $\A[\SU_q(2)]$ we have \begin{equation}\label{su gens}\mathbf{x}=\begin{pmatrix}a&b\\c&d\end{pmatrix}=\begin{pmatrix}a & -qc^* \\
c & a^*\end{pmatrix}.\end{equation} The algebra relations become
\begin{equation}\label{newrels}ac=qca,
~~ac^*=qc^*a, ~~cc^*=c^*c, ~~aa^* + q^2 cc^*=1, ~~a^*a + c^*
c=1,\end{equation} together with their conjugates.
On generators, the counit is
$\ep(a)=\ep(a^*)=1$, $\ep(c)=\ep(c^*)=0$ and the antipode is now
$S(a)=a^*$, $S(a^*)=a$, $S(c)=-qc$, $S(c^*)=-\iq c^*$, while the coproduct now reads
$\Delta(a)=a\otimes a -q c^*\otimes c$,  $\Delta(c)=c \otimes a + a^* \otimes c$ and
$\Delta(a^*)=a^*\otimes a^* -q c\otimes c^*$,  $\Delta(c^*)=c^* \otimes a^* + a \otimes c^*$.

\subsection{The quantum universal enveloping algebra
$\usu$}\label{se:envalg} The quantum universal enveloping algebra
$\usu$ is the unital $*$-algebra generated by the four elements $K$,
$K^{-1}$, $E$, $F$, with $KK^{-1}=K^{-1}K=1$, subject to the
relations
\begin{equation}\label{surels}K^{\pm 1} E=q^{\pm 1} EK^{\pm 1},\qquad
K^{\pm 1} F=q^{\mp 1} FK^{\pm 1},\qquad
[E,F]=(q-\iq)^{-1}\left(K^2-K^{-2}\right)\end{equation} and the
$*$-structure $$K^*=K,\qquad E^*=F,\qquad F^*=E.$$ It becomes a Hopf
$*$-algebra when equipped with the coproduct $\Delta$ and counit
$\ep$ defined on generators by $$\Delta(K^{\pm 1})=K^{\pm 1}\otimes
K^{\pm 1},\quad \Delta(E)=E\otimes K+K^{-1}\otimes E,\quad
\Delta(F)=F\otimes K+K^{-1}\otimes F,$$ $$\ep(K)=1,\qquad
\ep(E)=0,\qquad \ep(F)=0,$$ and with antipode $S$ defined by
$S(K)=K^{-1}$, $S(E)=-qE$, $S(F)=-\iq F$ on generators. The maps
$\Delta$, $\ep$ are extended as $*$-algebra maps, whereas $S$
extends as a $*$-anti-algebra map. From the relations
\eqref{surels}, one finds that the quadratic Casimir element
\begin{equation}\label{casimir}
C_q:=FE+(q-\iq)^{-2}\left(qK^2 - 2 +\iq K^{-2}\right) - \tfrac{1}{4}
\end{equation}
generates the centre of the algebra $\usu$.

The finite-dimensional irreducible $*$-representations $\pi_j$ of
$\usu$ are indexed by a half-integer $j=0,1/2,1,3/2,\ldots$ called
the {\em spin} of the representation. Explicitly, these
representations are given by
\begin{align}
\label{irrep1}\pi_j(K) |j,m\ra &=q^m |j,m\ra, \\
\pi_j(F) |j,m\ra&=\left([j-m][j+m+1]\right)^{1/2} |j,m+1\ra, \nn \\
\pi_j(E) |j,m\ra&=\left([j-m+1][j+m]\right)^{1/2} \nn |j,m-1\ra,
\end{align}
where the vectors $|j,m\ra$ for $m=-j,-j+1,\ldots,j-1,j$ form an
orthonormal basis of the $(2j+1)$-dimensional irreducible
$\usu$-module $V^j$. Moreover, $\pi_j$ is a $*$-representation with
respect to the Hermitian inner product on $V^j$ for which the
vectors $|j,m\ra$ are orthonormal. In each representation, the
Casimir $C_q$ of~\eqref{casimir} acts as a multiple of the
identity, with constant given by
\begin{equation}\label{cas-spec}
\pi_j(C_q)=[j+\tfrac{1}{2}]^2 - \tfrac{1}{4}
\end{equation}
as one may easily verify by direct computation.

The Hopf $*$-algebras $\A(\SU_q(2))$ and $\usu$ are dually paired
{\em via} a bilinear pairing
\begin{equation}\label{pair}
(\,\cdot\,,\,\cdot\,):\,\usu\times \A[\SU_q(2)]\to \C
\end{equation}
which is non-degenerate.
It is defined on generators by
\begin{equation*}(K,a)=q^{-1/2},\qquad
(K^{-1},a)=q^{1/2},\qquad (K,d)=q^{1/2},\qquad
(K^{-1},d)=q^{-1/2},\end{equation*} \begin{equation*} (E,c)=1,\qquad
(F,b)=1,\end{equation*}
with all other combinations of generators
pairing to give zero. The pairing is extended to products of
generators {\em via} the requirements
\begin{equation}\label{extend}(\Delta(X),p_1\otimes p_2)=(X,p_1p_2),\qquad
(X_1X_2,p)=(X_1\otimes X_2,\Delta(p)),\end{equation}
\begin{equation*}(X,1)=\ep(X),\qquad (1,p)=\ep(p)\end{equation*} for
all $X,X_1,X_2\in \usu$ and all $p,p_1,p_2\in\A[\SU_q(2)]$. It is
compatible with the antipode and the $*$-structures in the sense
that, for all $X\in\usu$,
$p\in \A[\SU_q(2)]$,
\begin{equation}
(S(X),p)=(X,S(p)),\qquad
(X^*,p)=\overline{(X,(S(p))^*)},\qquad
(X,p^*)=\overline{((S(X))^*,p)} .
\end{equation}

Using the pairing, there is a canonical left
action of $\usu$ on $\ASU$ defined by
\begin{equation}\label{leftact} \tr:\usu\times\ASU\to \ASU,\qquad
X\tr p:=p\o (X,p\t)\end{equation}
where $X\in\usu$, $p\in\ASU$ and
$\Delta(p)=p\o\otimes p\t$ denotes the coproduct on $\ASU$. In
particular, this action works out on generators to be
\begin{eqnarray}\label{leftgens1}
E\tr a=b,\quad E\tr c=d,\quad F\tr b=a,\quad F\tr d=c,  \\
K^{\pm 1}a=q^{\pm1/2}a,\quad K^{\pm 1}c=q^{\pm 1/2}c,\quad
K^{\pm 1}b=q^{\mp 1/2}b,\quad K^{\pm 1}d=q^{\mp 1/2}d, \nn \\
E\tr b=0,\quad E\tr d=0,\quad F\tr a=0,\quad F\tr c=0 \nn .
\end{eqnarray}
This action makes $\ASU$ into a left $\usu$-module $*$-algebra, in the sense that
$$X\tr(p_1p_2)=(X\o\tr p_1)(X\t\tr p_2),\qquad X\tr 1=1,\qquad X\tr
p^*=((S(X))^*\tr p)^*$$ for all $p,p_1,p_2\in\ASU$, $X\in \usu$.
There is also a canonical right action of $\usu$ on $\ASU$, defined
by
\begin{equation}\label{rightact}
\tl: \ASU \times \usu \to \ASU,\qquad
p \tl X := (X,p\o) p\t
\end{equation}
for $X\in\usu$ and $p\in\ASU$, with properties similar to those for
the left action. These two canonical actions commute amongst one
another.

\subsection{Line bundles on the quantum sphere $S^2_q$}
The coordinate algebra $H:=\A[\U(1)]$ of the group $\U(1)$ is the
commutative unital $*$-algebra generated by $t,t^*$, subject to the
relations $tt^*=t^*t=1$. It is a Hopf algebra when equipped with
the coproduct, counit and antipode $$\Delta(t)=t\otimes t,\qquad
\ep(t)=1,\qquad S(t)=t^*,$$ extended as $*$-algebra maps.
There is a canonical Hopf algebra projection given on
generators by
\begin{equation}\label{hopfproj}\pi:\A[\SU_q(2)]\to\A[\U(1)],\qquad \pi\begin{pmatrix}a & b \\
c & d\end{pmatrix}:=\begin{pmatrix}t & 0 \\
0 & t^*\end{pmatrix}.\end{equation}
Using this projection
a right coaction of $H=\A[\U(1)]$ on $P:=\A[\SU_q(2)]$ is defined by
\begin{equation}\label{r-coact}\delta_R:\A[\SU_q(2)]\to\A[\SU_q(2)]\otimes\A[\U(1)],\qquad
\delta_R(\sx_i{}^j):=\sx_i{}^\mu\otimes
\pi(\sx_\mu{}^j).\end{equation} In fact this coaction is the same
thing as a $\ZZ$-grading on $\A[\SU_q(2)]$ for which the generators
have degrees
\begin{equation}\label{grad}\textup{deg}(a)=\textup{deg}(c)=1,\qquad
\textup{deg}(b)=\textup{deg}(d)=-1.\end{equation} The subalgebra of
coinvariants under this coaction is denoted $\A[S^2_q]$,
$$\A[S^2_q]:=\{m\in\A[\SU_q(2)]~|~\delta_R(m)=m\otimes 1\}.$$ We shall frequently write $M:=\A[S^2_q]$. This algebra is precisely the subalgebra generated by elements of degree zero: it is the unital $*$-algebra generated by the elements
\begin{equation}\label{sphgens}b_+:=cd,\qquad
b_-:=ab,\qquad b_0:=bc\end{equation} subject to the relations
$$b_0b_\pm=q^{\pm 2}b_\pm b_0,\qquad
q^{-2}b_-b_+=q^2b_+b_-+(1-q^2)b_0,$$
$$b_+b_-=b_0(1+\iq b_0)$$ inherited from those
of $\A[\SU_q(2)]$. In the classical limit $q\rightarrow 1$, the
first line of relations becomes the statement that the algebra is
commutative, whereas the second line becomes the sphere relation for
the classical two-sphere $S^2$. The quantum sphere $S^2_q$ is
precisely the standard Podle\'{s} sphere of \cite{pod1}. The
canonical algebra inclusion $M\hookrightarrow P$ is well-known to be
a Hopf-Galois extension \cite{brmaj} and hence a quantum principal
bundle with universal differential calculi whose typical fibre is
determined by $H:=\A[\U(1)]$.

The coaction \eqref{r-coact} of $H$ on $\ASU$ is also used to define
a family of line bundles over the quantum sphere $S^2_q$, indexed by
$n\in\ZZ$:
$$
\mathcal{L}_n:=\left\{x\in \ASU~|~\delta_R(x)=x\otimes
t^{-n}\right\}.
$$
One has the decomposition \cite{mmnnu}
$$
\ASU=\bigoplus_{n\in\ZZ} \mathcal{L}_n .
$$
In particular $\mathcal{L}_0=\A[S^2_q]$ and one finds that $\mathcal{L}_n^*\cong
\mathcal{L}_{-n}$ and $\mathcal{L}_n\otimes_{\A[S^2_q]}\mathcal{L}_m\cong\mathcal{L}_{n+m}$
for each $n,m\in\ZZ$. Moreover,
$$
E\tr \mathcal{L}_n\subset \mathcal{L}_{n+2}, \qquad F\tr
\mathcal{L}_n\subset\mathcal{L}_{n-2},\qquad K^{\pm 1} \tr
\mathcal{L}_n\subset \mathcal{L}_n$$ for all $n\in \ZZ$, as can be
checked directly using~\eqref{leftgens1} and~\eqref{extend}.

It is known that each $\mathcal{L}_n$ is a finitely generated
projective (say) left $\A[S^2_q]$-module of rank one \cite{sw}. In
this way, we think of the module $\mathcal{L}_n$ as the space of
sections of a line bundle over $S^2_q$ with winding number $-n$.

\section{Differential Structure of the Quantum Hopf
Fibration}\label{se:diffs}

In this section we equip the quantum group $\SU_q(2)$ with a
four-dimensional bicovariant differential calculus, originally
described in \cite{wor}. Using this, the base space $S^2_q$ of the
Hopf fibration inherits a three-dimensional differential calculus
which was originally described in \cite{pod2}, although we describe
it here in terms which are more compatible with the principal bundle
structure. Finally we show that $S^2_q$ is a framed quantum
manifold, in the sense that its cotangent bundle is a vector bundle
associated to the Hopf fibration $\SU_q(2)\to S^2_q$.

\subsection{Differential structure on $\SU_q(2)$}In the following we write $\ep_P$ for the counit of the Hopf algebra
$P:=\A[\SU_q(2)]$. In terms of the matrix elements in~\eqref{su
gens}, we define $I_P$ to be the right ideal of
$P^+:=\textup{Ker}\,\ep_P$ generated by the nine elements
\begin{equation}\label{idgens}b^2,~c{}^2,~b(a-d),~c(a-d),~a^2+q^{2}d^2-(1+q^2)(ad+\iq
bc),\end{equation}
$$\st b,~\st c,~\st(a-d),~\st(q^2a+d-(q^2+1)),$$
where $\st:=q^2a+d-(q^3+q^{-1})$. As discussed in \S\ref{diffs},
this ideal defines a left-covariant first order differential
calculus on $\SU_q(2)$, which we denote by $\Omega^1P$. In fact, one
checks that $I_P$ is stable under the right adjoint coaction
$\textup{Ad}_R$ and so this calculus is bicovariant under left and
right coactions of $\A[\SU_q(2)]$. It is precisely the 4D${}_+$
calculus on $\SU_q(2)$ introduced in \cite{wor}: indeed, one may
check that the space $\Lambda^1\cong P^+/I_P$ of left-invariant
one-forms is a four-dimensional vector space.

Following
\cite{ks:qgr}, we define elements $L_-$, $L_0$, $L_+$, $L_z$ of
$\usu$ by
\begin{equation*}
L_-:=q^{1/2}FK^{-1},\qquad L_+:=q^{-1/2}EK^{-1},
\end{equation*}
\begin{equation*}
L_0:=K^{2}+\nu^2\iq FE-1, \qquad
L_z:=K^{-2}-1.
\end{equation*} The vectors $L_0$ and $L_z$ are related
to the quantum Casimir \eqref{casimir} by
\begin{equation}\label{LCas}
(q-\iq)^2 \left( C_q + \tfrac{1}{4} - [\tfrac{1}{2}]^2 \right) =qL_0+\iq L_z .
\end{equation}
The elements $L_-$, $L_0$, $L_+$, $L_z$ act
upon $\ASU$ {\em via} the formula \eqref{leftact} and together
provide a basis for the quantum tangent space $\mathcal{T}_P$ of the calculus. Note
in particular that the element $C_q-\ep_P(C_q)1$ is also an element
of $\mathcal{T}_P$.

Let $\{\omega_-, \omega_0, \omega_+, \omega_z\}$ be a basis of the
space of left-invariant one-forms $\Lambda^1$ such that
$(L_j,\omega_k)=\delta_{jk}$ for $j,k=-,0,+,z$. As given in
\cite{ks:diff}, the bimodule relations in the calculus $\Omega^1P$
with respect to these one-forms are:
\begin{eqnarray}\label{bimods1}
\omega_-\begin{pmatrix}a&b\\c&d\end{pmatrix}=\begin{pmatrix}a&b\\c&d\end{pmatrix}\omega_-+\nu^2\iq\begin{pmatrix}b&0\\d&0\end{pmatrix}\omega_0;  \\
\omega_+\begin{pmatrix}a&b\\c&d\end{pmatrix}=\begin{pmatrix}a&b\\c&d\end{pmatrix}\omega_++\nu^2\iq\begin{pmatrix}0&a\\0&c\end{pmatrix}\omega_0;\qquad
\omega_0\begin{pmatrix}a&b\\c&d\end{pmatrix}=\begin{pmatrix}\iq
a&qb\\ \iq c&qd\end{pmatrix}\omega_0; \nn \\
\omega_z\begin{pmatrix}a&b\\c&d\end{pmatrix}=\begin{pmatrix}0&a\\0&c\end{pmatrix}\omega_-+\nu^2\iq\begin{pmatrix}a&0\\c&0\end{pmatrix}\omega_0+
\begin{pmatrix}b&0\\d&0\end{pmatrix}\omega_+
+\begin{pmatrix}qa&\iq b\\qc&\iq
d\end{pmatrix}\omega_z. \nn
\end{eqnarray}
In these terms, the exterior
derivative $\D:\A[\SU_q(2)]\to\Omega^1P$ has the form
\begin{equation}\label{ext}\D p=(L_-\tr p)\omega_- +(L_0\tr
p)\omega_0+(L_+\tr p)\omega_+ + (L_z\tr p)\omega_z,\qquad
p\in\ASU,\end{equation} where $\tr$ is the left action of $\usu$ on
$\A[\SU_q(2)]$ defined in~\eqref{leftact}. By using the
formul{\ae} \eqref{leftgens1} to compute the
action of $L_0$, $L_z$, $L_+$, $L_-$ on the generators of $\ASU$ and
then substituting into~\eqref{ext}, one obtains the explicit
expressions
\begin{align}\label{diffs1}
\D a&=\left( \iq-1+\nu^2\iq\right) a\omega_0+b\omega_++(q-1)a\omega_z,
\\ \D b&=a\omega_-+(q-1)b\omega_0+(\iq-1)b\omega_z, \nn \\
\D c&=\left(\iq-1+\nu^2\iq\right)c\omega_0+d\omega_++(q-1)c\omega_z, \nn \\
\D d&=c\omega_-+(q-1)d\omega_0+(\iq-1)d\omega_z \nn
\end{align} for the
differentials of the matrix generators of $\ASU$ in terms of these
left-invariant one-forms.

\subsection{Framed manifold structure of $S^2_q$} Next we use
\S\ref{se:frames} to compute the cotangent bundle $\Omega^1S^2_q$ of
the base space $S^2_q$ of the Hopf fibration as an associated vector
bundle. As before, we write $P=\ASU$ for the algebra of functions on
the total space of the Hopf fibration, $M=\A[S^2_q]$ for the algebra
of functions on the base and $H=\A[\U(1)]$ for the structure quantum
group. Recall the right coaction $\delta_R:P\to P\otimes H$ defined
in~\eqref{r-coact} and the canonical projection $\pi:P\to H$
defined in~\eqref{hopfproj}.

The differential calculus on $P$ is taken to be the four-dimensional
bicovariant calculus $\Omega^1P$ defined in the previous section;
it is defined in terms of the $\textup{Ad}_R$-invariant
ideal $I_P$ generated by the elements in~\eqref{idgens}. Now
writing $\ep_H$ for the counit of $H$, we obtain a bicovariant
differential calculus $\Omega^1H$ on $H=\A[\U(1)]$ by projecting the
ideal $I_P$ to obtain an ideal $I_H:=\pi(I_P)$ of
$\textup{Ker}\,\ep_H$. As such, $I_H$ is generated by the three
elements
\begin{equation}\label{onegens}t^2+q^{2}t^*{}^2-(1+q^2),\quad\st(t-t^*),\quad\st(q^2t+t^*-(q^2+1)),\end{equation}
again with $\st=q^2t+t^*-(q^3+\iq)$, where $t,t^*$ are the
generators of $H$.

\begin{lem}\label{le:fibcal}The calculus $\Omega^1H$ is
one-dimensional. It is spanned as a left module by the
left-invariant one-form $\omega_t:=t^*\D t$ and has bimodule
relations
$$\omega_t t=qt\omega_t, \qquad \omega_t t^*=\iq t^*\omega_t,$$
where $t,t^*$ are the generators of $H=\A[\U(1)]$.\end{lem}

\proof We define an equivalence relation $\sim$ on $H^+$ by $x\sim
y$ if and only if $x-y\in I_H$. By taking a linear combination of
the generators in \eqref{onegens}, one finds in particular that
$(t-1)+q(t^*-1)\sim 0$, which is our key equivalence. Using it, one
deduces that
\begin{align*}
t^2&=(t+1)(t-1)+1\sim -q(t+1)(t^*-1)+1 =-q(t^*-t)+1\sim (q+1)(t-1)+1,\\
t^*{}^2&=(t^*+1)(t^*-1)+1\sim -\iq(t^*+1)(t-1)+1\\
&\qquad\qquad\qquad\qquad\qquad =-\iq(t-t^*)+1\sim - \iq(1+\iq)(t-1)+1,\end{align*}
so that every quadratic polynomial in $t$, $t^*$ and $1$ is
equivalent to a linear combination of $t-1$ and $t^*-1$. By
induction any polynomial in $t$ is equivalent
to such a linear combination. Applying the key equivalence once more
tells us that we can always eliminate $t^*-1$. Thus we take $t-1$ as
a representative of the quotient space $H^+/I_H$ and
$\omega_t:=r^{-1}(1\otimes(t-1))$ as the corresponding
left-invariant one-form, which spans the calculus $\Omega^1H$ as a
left $H$-module. To obtain the bimodule relations, we compute for
example that
$$
\omega_tt=((t^*-1)\otimes \lfloor t-1\rfloor)t=(1-t)\otimes
\lfloor t^2-t\rfloor=qt(t^*-1)\otimes \lfloor
t-1\rfloor=qt\omega_t,
$$
where $\lfloor ~ \rfloor$ denotes an
equivalence class modulo $I_H$. The first and last equalities use
the definition of the map $r$ and the middle equality uses the
bimodule structure \eqref{right module calculus}.\endproof

The differential calculus $\Omega^1M$ on the base of the fibration
is defined by restricting the calculus $\Omega^1P$ to $M$. This
means that it is defined as the quotient
$\Omega^1M:=\widetilde\Omega^1M/N_M$, where $N_M$ is the
$M$-$M$-bimodule $N_M:=N_P\cap\widetilde\Omega^1M$. We postpone the
computation of generators and relations for $\Omega^1M$ and observe
that for now we have the following expressions for the exterior
derivative on $M$ in terms of the left-invariant one-forms
$\omega_\pm, \omega_0$.

\begin{lem}\label{sphdiffs} The exterior derivative $\D$ acts on $M=\A[S^2_q]$ as
\begin{equation}\begin{pmatrix}
\D b_+\\\D b_0\\\D
b_-\end{pmatrix}=\begin{pmatrix}d^2&\mu\nu^2\iq cd&qc^2\\
db&\nu^2\iq\left( 1+\mu bc\right)&ac\\b^2&\mu\nu^2\iq
ab&qa^2
\end{pmatrix}\begin{pmatrix}\omega_+\\\omega_0\\\omega_-\end{pmatrix}
\end{equation}
in terms of the generators $b_\pm, b_0$ of $M$ given in~\eqref{sphgens}.
\end{lem}

\proof This follows from direct computation. For example, to compute
$\D b_+$ the Leibniz rule yields $$\D b_+= \D(cd)=(\D c)d + c(\D
d).$$ One uses the expressions \eqref{diffs1} to
rewrite $\D c$, $\D d$ in terms of $\omega_\pm$ and $\omega_0$, then
the bimodule relations in Eqs.~\eqref{bimods1}  to
collect all coefficients to the left. Combining together alike terms
yields the expression as stated. The same method works for
computing $\D b_0$ and $\D b_-$.\endproof

\begin{lem} With $P$, $H$ and $M$ as above, the differential calculi
$\Omega^1P$, $\Omega^1H$ and $\Omega^1M$ satisfy the compatibility
conditions of~\eqref{hom bundle}.\end{lem}

\proof The relation $\pi(I_P)=I_H$ holds by definition of the
calculus on $H$. It is sufficient to verify the
$\textrm{Ad}_R$-condition in \eqref{hom bundle} on generators: one
finds that \begin{align*}(\textup{id} \otimes
\pi)\textup{Ad}_R(c^2)&=c^2\otimes t^4,\qquad ~~(\textup{id} \otimes
\pi)\textup{Ad}_R(c(a-d))=c(a-d)\otimes
t^2,\\
(\textup{id} \otimes \pi)\textup{Ad}_R(b^2)&=b^2\otimes
t^*{}^4,\qquad (\textup{id} \otimes
\pi)\textup{Ad}_R(b(a-d))=b(a-d)\otimes t^*{}^2,\\
(\textup{id} \otimes \pi)\textup{Ad}_R(\st c)&=\st c\otimes
t^2,\qquad\quad (\textup{id} \otimes \pi)\textup{Ad}_R(\st b)=\st
b\otimes t^{-2},\end{align*} with all other generators coinvariant
under the map $(\textup{id} \otimes \pi)\textup{Ad}_R$.\endproof

This means that we may apply \S\ref{se:frames} to express $S^2_q$ as
a framed quantum manifold. The framing comodule $V$ is computed as
follows. Clearly $P^+\cap M$ is equal to $M^+=\textup{Ker}\,\ep_M$,
the restriction of the counit $\ep_P$ to the subalgebra $M$. In our
case, with $M=\A[S^2_q]$ being generated by $b_\pm, b_0$, we have
that $M^+=\la b_0,b_\pm\ra$ as a right ideal. To compute $I_P\cap M$
we note that, since the generators $b(a-d)$, $c(a-d)$,
$a^2+q^{2}d^2-(1+q^2)(ad+\iq bc)$, $\st b$, $\st c$, $\st(a-d)$,
$\st(q^2a+d-(q^2+1))$ are not of homogeneous degree, the ideal that
each of them generates has no intersection with $M$. 
Thus we concentrate on the generators $b^2$, $c^2$ of $I_P$. The
elements of degree zero in $\la b^2\ra$ include $b^2\{ a^2,ac,c^2\}$
and so we see that $b_-^2,b_-b_0,b_0^2$ all lie in $I_P\cap M$.
Similarly, from the ideal $\la c^2\ra$ we see that $b_+^2$ and
$b_+b_0$ are also in $I_P\cap M$. From this discussion we obtain
\begin{equation}\label{comod}V=\la b_0, b_\pm\ra/\la
b_\pm^2,b_0^2,b_\pm b_0\ra.\end{equation} Hence $V$ is
three-dimensional with representatives $b_\pm$ and $b_0$. We compute
the right coaction of $H$ on $V$ from~\eqref{frame} as
$$
\Delta_R(b_+)=cd\otimes S\pi(d^2)=b_+\otimes
t^2,\quad \Delta_R(b_-)=ab\otimes S\pi(a^2)=b_-\otimes t^*{}^2,
$$
$$
\Delta_R(b_0)=bc\otimes 1=b_0\otimes 1.
$$ Hence
$V=\C\oplus\C\oplus\C$ and the associated bundle
$$\E=\mathcal{L}_{-2}\oplus\mathcal{L}_0\oplus\mathcal{L}_{+2}=\A[\SU_q(2)]_2\oplus
\A[\SU_q(2)]_0\oplus\A[\SU_q(2)]_{-2}$$ is the direct sum of the
line bundles over $S^2_q$ with winding numbers $-2$, $0$ and $2$. This yields
the following theorem.

\begin{thm}\label{th:framing}The homogeneous space $S^2_q$ is a framed quantum
manifold with cotangent bundle
$$\Omega^1S^2_q\cong\mathcal{L}_{-2}\oplus\mathcal{L}_0\oplus\mathcal{L}_{+2}.$$
The isomorphism is given by the soldering form
\begin{align*}\theta(b_+)&=q^2c^2\D b_--q\mu ac\D
b_0+a^2\D b_+=\omega_+,\\ \theta(b_0)&=-qdc\D b_-+(1+\mu bc)\D b_0-\iq ba\D b_+=\nu^2\iq \omega_0,\\
\theta(b_-)&=d^2\D b_--\iq\mu bd\D b_0+q^{-2}b^2\D
b_+=q\omega_-\end{align*} and makes $\Omega^1S^2_q$ projective as a
left $\A[S^2_q]$-module.\end{thm}

\proof The only remaining part is to compute the soldering form
$\theta(b_\pm)$, $\theta(b_0)$. We find the left coaction on
$M=\A[S^2_q]$ inherited from the coproduct on $\A[\SU_q(2)]$ to be
\begin{align*}\Delta_L(b_+)&=\Delta_L(cd)=c^2\otimes
b_-+cd\otimes (1+\mu b_0)+d^2\otimes b_+,\\
\Delta_L(b_0)&=\Delta_L(bc)=ca\otimes b_-+1\otimes
b_0+bc\otimes(1+\mu b_0)+db\otimes b_+,\\
\Delta_L(b_-)&=\Delta_L(ab)=a^2\otimes b_-+ab\otimes (1+\mu
b_0)+b^2\otimes b_+.\end{align*} In fact these coproducts were
already used in computing $\Delta_R$ above. This time we apply the antipode $S$
to the first tensor factor to obtain
$$\theta(b_+)=S(b_+\o)\D (b_+\t)=q^2c^2\D b_--q\mu ac\D
b_0+a^2\D b_+,$$ similarly for $\theta(b_-)$ and $\theta(b_0)$. This
yields the middle expressions as stated. We then insert the
expressions from Lemma~\ref{sphdiffs} to obtain
$\{\omega_+,\nu^2\iq\omega_0,q\omega_-\}$ for the values of the map
$\theta$. According to \S\ref{se:frames}, the map $\theta:V\to
P\Omega^1M$ is well-defined on $V$. In order to get one-forms on
$\A[S^2_q]$, one must multiply $\theta(b_-)$ by an element of degree
$2$, $\theta(b_+)$ by an element of degree $-2$ and $\theta(b_0)$ by
an element of degree zero. Moreover, every one-form is obtained in
this way. This yields the isomorphism as stated. Since all line
bundles $\mathcal{L}_n$ are projective, so is
$\Omega^1S^2_q$.\endproof

The above also shows that the exterior derivative $\D$ in the calculus
$\Omega^1S^2_q$ is given by restriction of the expression in~\eqref{ext}, namely
\begin{equation}\label{res ext} \D m=(L_-\tr m)\omega_-+(L_0\tr
m)\omega_0+(L_+\tr m)\omega_+,\qquad m\in \A[S^2_q].
\end{equation}
We stress that $L_{\mp} \tr m \in {\mathcal L}_{\pm 2}$ rather then being element in $\A[S^2_q]$.
Of course, from~\eqref{ext} combined with the fact that the
vertical vector field $L_z$ obeys $L_zm=0$ for all $m\in \A[S^2_q]$,
we already expected this to be the case. From
Theorem~\ref{th:framing} we know that $\Omega^1S^2_q$ is spanned as
a left module by
\begin{align}\label{hol1}
\{d^2,db,b^2\}\, \omega_+&:=\{\p_+b_+,\p_+b_0,\p_+b_-\}, \\
\nu^2\iq\{\mu cd,1+\mu bc,\mu ab\}\, \omega_0&:=\{\p_0b_+,\p_0b_0,\p_0b_-\}, \nn \\
\{qc^2,ac,qa^2\}\, \omega_-&:=\{\p_-b_+,\p_-b_0,\p_-b_-\}. \nn
\end{align}
The bimodule relations in the calculus $\Omega^1S^2_q$ are in
general quite complicated to compute directly, but we can use the
expressions in Eqs.~\eqref{hol1} to break them into smaller pieces
which are much easier to work with.

\begin{cor}\label{co:birels} The cotangent bundle $\Omega^1S^2_q$
has first order differential sub-calculi $$\Omega^1_+\cong
\mathcal{L}_{-2}\oplus\mathcal{L}_0,\qquad \Omega^1_0\cong
\mathcal{L}_0, \qquad \Omega^1_-\cong
\mathcal{L}_0\oplus\mathcal{L}_{+2}$$ with differentials given by
$\D_+:=\p_++\p_0$, $\D_0:=\p_0$ and $\D_-:=\p_0+\p_-$ respectively.
These calculi obey the bimodule relations
\begin{align*}
\p_+b_+\left\{\begin{array}{l}b_+\\
b_0\\b_-\\\, \end{array}\right.&=\left\{
\begin{array}{l}q^{-2}b_+(\p_+b_+)+q^{-3}\mu^{-1}b_+(\p_0b_+)\\q^{-4}
b_0(\p_+b_+)+\mu^{-1}q^{-2}(1+q^{-3}b_0)(\p_0b_+)
\\ q^{-2}b_-(\p_+b_+)-(q^2-q^{-2})b_+(\p_+b_-)+\p_0 b_0\\ \qquad\qquad +(q^2-q^{-2})^{-1}\left( q^{-2}b_-(\p_0b_+)-b_+(\p_0b_-)\right)-\iq\nu b_+(\p_0b_-),\end{array}\right.
\\ ~\\
\p_+b_0\left\{\begin{array}{l}b_+\\
b_0\\b_-\end{array}\right.&=\left\{
\begin{array}{l}b_+(\p_+b_0)+q^{-3}\mu^{-1}b_0(\p_0b_+)\\q^{-2}
b_0(\p_+b_0)+q^{-2}\mu^{-1}b_+(\p_0b_-)
\\ q^{-2}b_-(\p_+b_0)-q^{-1}\nu b_0(\p_+b_-)+q^{-2}(1+q^{-1}b_0)(\p_0b_-),\end{array}\right.
\\ ~\\
\p_+b_-\left\{\begin{array}{l}b_+\\
b_0\\b_-\end{array}\right.&=\left\{
\begin{array}{l}q^2b_+(\p_+b_-)+(q^2-q^{-2})^{-1}\left(q^{2}b_-(\p_0b_+)-b_+(\p_0b_-)\right)\\
b_0(\p_+b_-)+\iq\mu^{-1}b_0(\p_0b_-)
\\ q^{-2}b_-(\p_+b_-)+q^{-3}\mu^{-1}b_-(\p_0b_-),\end{array}\right.
\\ ~\\
\p_-b_+\left\{\begin{array}{l}b_+\\
b_0\\b_-\end{array}\right.&=\left\{
\begin{array}{l}q^2b_+(\p_-b_+)+q^{3}\mu^{-1}b_+(\p_0b_+)\\
b_0(\p_-b_+)+q\mu^{-1}b_0(\p_0b_+)
\\ q^{-2}b_-(\p_-b_+)+(q^2-q^{-2})^{-1}\left(b_-(\p_0b_+)-q^2b_+(\p_0b_-)\right),\end{array}\right.
\\ ~\\
\p_-b_0\left\{\begin{array}{l}b_+\\
b_0\\b_-\end{array}\right.&=\left\{
\begin{array}{l}q^2b_+(\p_-b_0)+q\nu b_0(\p_-b_+)+\mu^{-1}(1+q b_0)(\p_0b_+)\\q^2
b_0(\p_-b_0)+\mu^{-1}b_-(\p_0b_+)\\ b_-(\p_-b_0)+q^3\mu^{-1}
b_0(\p_0b_-),\end{array}\right.
\\ ~\\
\p_-b_-\left\{\begin{array}{l}b_+\\ \, \\
b_0\\b_-\end{array}\right.&=\left\{
\begin{array}{l} q^2b_+(\p_-b_-)+(q^2-q^{-2})b_-(\p_-b_+)+q^2\p_0b_0 \\ \qquad\qquad + (q^2-q^{-2})^{-1}\left(b_-(\p_0b_+)-q^2b_+(\p_0b_-)\right)+q\nu b_-(\p_0 b_+) \\q^4 b_0(\p_-b_-)+\mu^{-1}(1+q^3
b_0)(\p_0b_-)
\\ q^2b_-(\p_-b_-)+q^3\mu^{-1}b_-(\p_0b_-).\end{array}\right.
\end{align*}\end{cor}

\proof Using the expressions in Eqs.~\eqref{hol1} the bimodule relations in $\Omega^1S^2_q$ are easily determined from straightforward but laborious computation along the following lines.
From the bimodule relations in Eqs.~\eqref{bimods1} one finds that
$$\omega_+ \left\{\begin{array}{l}b_+\\b_0\\b_-\end{array}\right.
=\left\{\begin{array}{l} b_+ \omega_++\nu^2\iq c^2\omega_0\\
b_0\omega_++\nu^2q^{-1} ca\omega_0\\
b_-\omega_++\nu^2\iq a^2\omega_0,\end{array}\right.\qquad \omega_-
\left\{\begin{array}{l}b_+\\b_0\\b_-\end{array}\right.
=\left\{\begin{array}{l} b_+\omega_-+\nu^2d^2\omega_0\\
b_0\omega_-+\nu^2db\omega_0\\
b_-\omega_-+\nu^2 b^2\omega_0,\end{array}\right.$$ with $\omega_0$
commuting with each of $b_\pm,b_0$. Combining these with the algebra
relations in $\A[\SU_q(2)]$ yields the bimodule relations as stated,
together with
\begin{align*}\p_0b_+\left\{\begin{array}{l}b_+\\
b_0\\b_-\end{array}\right.&=\left\{
\begin{array}{l}b_+(\p_0b_+)\\q^{-2}
b_0(\p_0b_+)\\ q^{-2}b_-(\p_0b_+)-q^{-2}b_-(\p_0b_+)+b_+(\p_0b_-),
\end{array}\right.
\\ ~\\
\p_0b_0\left\{\begin{array}{l}b_+\\
b_0\\b_-\end{array}\right.&=\left\{
\begin{array}{l}q^2b_+(\p_0b_0)-q\mu^{-1}\nu(\p_0b_+)\\b_0(\p_0b_0) \\
q^{-2}b_-(\p_0b_0)+\iq\mu^{-1}\nu(\p_0b_-),
\end{array}\right.
\\ ~\\
\p_0b_-\left\{\begin{array}{l}b_+\\
b_0\\b_-\end{array}\right.&=\left\{
\begin{array}{l}q^2b_+(\p_0b_-)+b_-(\p_0b_+)-q^{-2}b_+(\p_0b_-)\\q^2
b_0(\p_0b_-)\\ b_-(\p_0b_-).
\end{array}\right.
\end{align*} The fact that $\Omega^1_+=\mathcal{L}_{-2}\oplus \mathcal{L}_0$, $\Omega^1_0=\mathcal{L}_0$ and
$\Omega^1_-=\mathcal{L}_0\oplus\mathcal{L}_{+2}$ close as
sub-bimodules is now clear by inspection. The Leibniz rules for the
differentials $\D_+$, $\D_0$ and $\D_-$ follow from the Leibniz rule
for $\D$ and the direct sum decomposition of
$\Omega^1S^2_q$.\endproof

\begin{cor}The one-forms in the calculus $\Omega^1S^2_q$ enjoy the
relations
\begin{eqnarray*}\p_+b_0=q^{-2}
b_-(\p_+b_+)-q^2b_+(\p_+b_-),\qquad b_0b_-(\p_+b_+)=
q^3(1+qb_0)b_+(\p_+b_-), \\
\p_-b_0=b_+(\p_-b_-)-q^{-4}b_-(\p_-b_+),\qquad
b_0b_+(\p_-b_-)=q^{-3} (1+\iq b_0)b_-(\p_-b_+), \\
b_0\p_0b_0=-q\mu\nu^{-1} b_-(\p_0b_+)+\iq\mu\nu^{-1} b_+(\p_0b_-), \\
b_+(\p_0b_0)=(\mu^{-1}+q^{-2}b_0)\p_0b_+,\qquad b_-(\p_0b_0)=(\mu^{-1}+q^{2}b_0)\p_0b_+.
\end{eqnarray*}
\end{cor}

\proof These are obtained in analogy with the proof of
Cor.~\ref{co:birels}, from the relations in $\A[\SU_q(2)]$
acting on $\omega_\pm$ and $\omega_0$. One finds the relations as
stated, together with
\begin{align*}
b_+(\p_+b_-)&=q^{-1}b_0(\p_+b_0),\qquad
b_-(\p_+b_+)=q^2(1+qb_0)(\p_+b_0),\\
b_-(\p_-b_+)&=q^2b_0(\p_-b_0),\qquad~\; b_+(\p_-b_-)=\iq(1+\iq
b_0)(\p_-b_0).
\end{align*} There are other relations involving the
differential $\p_0$, but they are quite complicated (since the
sphere relation in $\A[S^2_q]$ does not explicitly involve the unit)
and are not particularly illuminating, so we shall not give them
here.\endproof

Finally, we use Theorem~\ref{th:framing} to compute the
differentials $\p_\pm$ and $\p_0$ in terms of the exterior
derivative $\D$. Using the algebra relations in $\A[\SU_q(2)]$ and
the expressions in Eqs.~\eqref{hol1} we find that
\begin{align*}\p_+b_+&=\iq b_+^2\D b_--\mu b_+(1+\iq b_0)\D
b_0+(1+\iq b_0)^2\D b_++q^{-2}\nu b_+b_-\D b_+, \\
\p_+b_0&=qb_+b_0\D
b_--\mu b_+b_-\D b_0+q^{-2}(1+\iq b_0)b_-\D b_+,\\
\p_+b_-&=q^2b_0^2\D b_--\iq\mu b_-b_0\D b_0+q^{-3}b_-^2\D b_+,\\
\p_0 b_+&=-\mu b_+^2\D b_-+\mu b_+(1+\mu b_0)\D b_0-q^{-2}\mu
b_+b_-\D b_+,\\ \p_0 b_0&=(1+\mu b_0)\left( -b_+\D b_- +(1+\mu
b_0)\D b_0-q^{-2} b_-\D b_+\right),\\ \p_0 b_-&=-\mu b_-b_+\D
b_-+\mu b_-(1+\mu b_0)\D b_0-q^{-2}\mu b_-^2\D b_+,\\
\p_-b_+&=qb_+^2\D b_--\iq \mu b_0b_+\D b_0+q^{-2}b_0^2\D b_+,
\\ \p_- b_0&=(1+qb_0)b_+\D b_--q\mu b_0(1+qb_0)\D b_0+q^{-2}b_-b_0\D
b_+,\\ \p_-b_-&=\left((1+qb_0)^2+\nu b_-b_+\right)\D b_--\mu b_-(1+q
b_0)\D b_0+\iq b_-^2\D b_+.\end{align*} These expressions may now be
used to compute the full bimodule structure of the calculus
$\Omega^1S^2_q$ in terms of the differential $\D$, as well as the
deeper structure of the noncommutative Riemannian geometry of this
calculus, along similar lines to \cite{Ma:sphere}. However, since
our objective is to study the spin geometry of the calculus, we have
all we need and so we shall not pursue these directions here.

\section{The Spectral Geometry of $S^2_q$}

\label{se:spectral} In this section we give the `three-dimensional'
differential calculus $\Omega^1S^2_q$ by a spectral triple on
$S^2_q$. This means equipping $S^2_q$ with a spinor bundle
$\mathcal{S}$ and a Dirac operator $D$ which together implement the
exterior derivative $\D$ for $\Omega^1S^2_q$. We then equip this
spectral triple with a real structure for which the commutant
property and the first order condition for the Dirac operator are
satisfied up to infinitesimals of arbitrary order, in parallel with
the results of \cite{dalw} for the `two-dimensional' calculus on
$S^2_q$.

\subsection{Background on spectral triples}\label{se:spectrips}
We recall  briefly the notion of a spectral triple~\cite{ac:book}.

\begin{defn} A {\em unital spectral triple} $(A,\mathcal{H},D)$ consists of a
complex unital $*$-algebra $A$, faithfully $*$-represented
by bounded operators on a (separable) Hilbert space $\mathcal{H}$,
and a self-adjoint operator $D:\mathcal{H}\to\mathcal{H}$ (the Dirac
operator) with the following properties:
\begin{enumerate}[\hspace{0.5cm}(i)] \item the resolvent
$(D-\lambda)^{-1}$, $\lambda \notin \RR$, is a compact operator on
$\mathcal{H}$; \item for all $a\in A$ the commutator $[D,\pi(a)]$
is a bounded operator on $\mathcal{H}$.\end{enumerate}
A spectral triple $(A,\mathcal{H},D)$ is called {\em even} if there
exists a $\mathbb{Z}_2$-grading of $\mathcal{H}$, {\em i.e.} an
operator $\Gamma:\mathcal{H}\to\mathcal{H}$ with $\Gamma=\Gamma^*$
and $\Gamma^2=1$, such that $\Gamma D+D\Gamma=0$ and $\Gamma
a=a\Gamma$ for all $a\in A$. Otherwise the spectral triple is said
to be {\em odd}.
\end{defn}

With $0<n<\infty$, the Dirac operator $D$ is said to be {\em
$n^+$-summable} if $(D^2+1)^{-1/2}$ is in the Dixmier ideal
$\mathcal{L}^{n^+}(\mathcal{H})$.
The {\em metric dimension} of the spectral triple
$(A,\mathcal{H},D)$ is defined to be the infimum of the set of all
$n$, such that $D$ is $n^+$-summable.

Given a spectral triple $(A,\mathcal{H},D)$, one associates to it a
canonical first order differential calculus $(\Omega^1_DA,\D_D)$.
In particular, the $A$-$A$-bimodule
$\Omega^1_DA$ is defined to be
\begin{equation}\label{commrep}
\Omega^1_DA:=\{\omega=\sum_j a^j_0[D,a^j_1]~|~a^j_0,a^j_1\in A\},
\end{equation}
with the differential $\D_D$ given by $\D_D a=[D,a]$ for $a\in A.$

The original definition \cite{ac:fncg} of a real structure on a spectral triple
$(A,\mathcal{H},D)$ was given by an anti-unitary operator
$J:\mathcal{H}\to\mathcal{H}$ with the properties $J^2=\pm 1$,
$JD=\pm DJ$ and
\begin{equation}\label{real}
[\pi(a),J\pi(b)J^{-1}]=0,\quad [[D,\pi(a)],J\pi(b)J^{-1}]=0,\qquad
a,b\in A.
\end{equation}
These are called the {\em commutant property} and the {\em first
order condition} respectively.

However, in many examples involving quantum spaces, one needs to
modify these conditions in order to obtain non-trivial spin
geometries \cite{dlps,dlss,dalw,ddl:iso}. Following the approach
there, we impose the weaker assumption that \eqref{real} holds only
up to infinitesimals of arbitrary order ({\em i.e.} up to compact
operators $T$ with the property that the singular values $s_k(T)$
satisfy $\textup{lim}_{k\to\infty}k^ps_k(T)=0$ for all $p>0$).

\begin{defn}\label{def:real} A {\em real structure} on a spectral triple
$(A,\mathcal{H},D)$ is an anti-unitary operator
$J:\mathcal{H}\to\mathcal{H}$ such that
$$J^2=\pm 1,\qquad JD=\pm DJ,
$$
\begin{equation}\label{infsmod}
[\pi(a),J\pi(b)J^{-1}]\in\mathscr{I},\quad
[[D,\pi(a)],J\pi(b)J^{-1}]\in\mathscr{I},\qquad a,b,\in A,
\end{equation}
where $\mathscr{I}$ is an operator ideal of infinitesimals of
arbitrary order. We say that the datum $(A,\mathcal{H},D,J)$ is a
{\em real spectral triple} (up to infinitesimals). If
$(A,\mathcal{H},D,\Gamma)$ is even and $J\Gamma=\pm\Gamma J$, we
call the datum $(A,\mathcal{H},D,\Gamma,J)$ an {\em even real
spectral triple} (up to infinitesimals).\end{defn}

The signs above depend on the so-called {\em KO-dimension} of the
triple. We shall only need the case where the KO-dimension is two;
then $J^2=-1$, $JD=DJ$ and $J\Gamma=-\Gamma J$.

\subsection{A Dirac operator on $S^2_q$}\label{se:spec} In order to
define a spectral triple on $S^2_q$, we need a spinor bundle over
$S^2_q$ and an associated Dirac operator, which we require should
recover the differential calculus $\Omega^1S^2_q$ {\em via} the
commutator representation defined in~\eqref{commrep}. Since the
differential calculus $\Omega^1S^2_q$ constructed in
Theorem~\ref{th:framing} is equivariant under a left coaction of
$\A[\SU_q(2)]$ and hence a right action of $\usu$, we are led to
consider spinor bundles and Dirac operators which are right
$\usu$-equivariant.

Guided by this principle, as well as by the spin
structure of the classical two-sphere $S^2$, for the $\A[S^2_q]$-module of spinors we take
$$
\mathcal{S}=\mathcal{S}_+\oplus\mathcal{S}_-:=\mathcal{L}_{-1}\oplus\mathcal{L}_{+1}.
$$
As right $\usu$-modules, the vector spaces $\mathcal{S}_\pm$ are
both isomorphic to the direct sum
\begin{equation}\label{spinrep}
V:=\bigoplus_{j\in\mathbb{N}+\tfrac{1}{2}}V^j
\end{equation}
over all irreducible $\usu$-modules $V^j$ with spin
$j\in\mathbb{N}+\tfrac{1}{2}$ a half-odd integer. A corresponding
basis for $V$ is then given by
$$
\{|j,m\ra~|~j\in\mathbb{N}+\tfrac{1}{2},m=-j,\ldots,j\} ,
$$
where the vectors $|j,m\ra$ span the irreducible $\usu$-module $V^j$
in Eqs.~\eqref{irrep1}. We denote the orthonormal bases of the two
different copies $\mathcal{S}_\pm$ of $V$ respectively by
\begin{equation}\label{basis}|j,m\ra_\pm,\qquad
j\in\mathbb{N}+\tfrac{1}{2},~m=-j,\ldots,j.
\end{equation}
We equip $\mathcal{S}$ with the inner product which makes this basis
orthonormal and write $\mathcal{H}$ for the corresponding Hilbert space completion
of $\mathcal{S}$.

As $\A[S^2_q]$-modules, the vector spaces $\mathcal{S}_\pm$ each
carry one of two inequivalent $\usu$-equivariant representations of
$\A[S^2_q]$,
$$\pi_{\pm}:\A[S^2_q]\to\textrm{End}(\mathcal{S}_\pm).$$ Recall that $\mathcal{S}_\pm$ are just the subspaces of $\ASU$ with
overall degrees $\mp 1$ with respect to the $\mathbb{Z}$-grading
\eqref{grad}, so the representations $\pi_{\pm}$ on
$\mathcal{S}_\pm$ are simply given by restricting the multiplication
in $\ASU$ to the appropriate degrees. However, it is possible to
describe these representations explicitly in terms of the basis
\eqref{basis} in the following way.

Indeed, the $\usu$-equivariant representations of $\A[S^2_q]$ on $V$
were already described in \cite{dalw,sw}. To be able to simply quote
them we make a change of generators, now writing
\begin{equation}\label{newgens}
x_1=- q^{1/2} \mu \, b_+,\qquad x_0-1=\mu \, b_0,\qquad
x_{-1}=-q^{-3/2}\mu \, b_-,
\end{equation}
where $b_\pm$, $b_0$ are the generators of $\A[S^2_q]$ defined in~\eqref{sphgens}, and $\mu=q+q^{-1}$.
With respect to these new generators, the algebra
relations of $\A[S^2_q]$ now read
\begin{align*}
x_{-1}(x_0-1)=q^{2}(x_0-1)x_{-1},\qquad x_1(x_0-1)=q^{-2}(x_0-1)x_1,\\
(q^{2}x_0+1)(x_0-1)=(q+\iq)x_{-1}x_1, \qquad
(q^{-2}x_0+1)(x_0-1)=(q+\iq)x_1x_{-1}.
\end{align*}
Then, with $N=\pm 1/2$, the two representations $\pi_\pm=\pi_{\pm
1/2}$ of $\A[S^2_q]$ on $\mathcal{S}_{\pm}$ have the form
\begin{multline}\label{irreps}
\pi_{N}(x_i)|j,m\ra_\pm =\alpha^-_i(j,m;N)|j-1,m+i\ra_\pm \\
\qquad\qquad+\alpha^0_i(j,m;N)|j,m+i\ra_\pm+\alpha^+_i(j,m;N)|j+1,m+i\ra_\pm,
\end{multline}
where the coefficients are determined by
\begin{align*}
\alpha^+_1(j,m;N)&=q^{-j+m}\left(\frac{[j+m+1][j+m+2]}{[2j+1][2j+2]}\right)^{1/2}\alpha_{N}(j+1), \\
\alpha^0_1(j,m;N)&=-q^{m+2}\left([2][j-m][j+m+1]\right)^{1/2}[2j]^{-1}\beta_N(j),\\
\alpha^-_1(j,m;N)&=-q^{j+m+1}\left(\frac{[j-m-1][j-m]}{[2j-1][2j]}\right)^{1/2}\alpha_N(j),\\
\alpha^+_0(j,m;N)&=q^{m}\left(\frac{[2][j-m+1][j+m+1]}{[2j+1][2j+2]}\right)^{1/2}\alpha_N(j+1),\\
\alpha^0_0(j,m;N)&=[2j]^{-1}\left([j-m+1][j+m]-q^{-2}[j-m][j+m+1]\right)\beta_N(j),\\
\alpha^-_0(j,m;N)&=q^{m}\left(\frac{[2][j-m][j+m]}{[2j-1][2j]}\right)^{1/2}\alpha_N(j),\\
\alpha^+_{-1}(j,m;N)&=q^{j+m}\left(\frac{[j-m+1][j-m+2]}{[2j+1][2j+2]}\right)^{1/2}\alpha_N(j+1),\\
\alpha^0_{-1}(j,m;N)&=q^{m}\left([2][j-m+1][j+m]\right)^{1/2}[2j]^{-1}\beta_N(j),\\
\alpha^-_{-1}(j,m;N)&=-q^{-j+m-1}\left(\frac{[j+m-1][j+m]}{[2j-1][2j]}\right)^{1/2}\alpha_N(j)
\end{align*}
(with the convention that
$\alpha^-_i(\tfrac{1}{2},\pm\tfrac{1}{2};N)=0$) and the real
numbers $\alpha_N(j)$, $\beta_N(j)$ are
\begin{align*}\alpha_N(j)&= \left( [2j+1][2j] \right)^{-1/2} \left( [2][j+N][j-N])^{1/2}  ([2j+1][2j] \right)^{1/2} \, q^{N},\\
\\~
\beta_N(j)&=
q^{-1} [2j+2]^{-1} \,
\left( \varepsilon q^{-\varepsilon}-(q-\iq)([j][j+1]-[\tfrac{1}{2}][\tfrac{3}{2}] \right)  ,
\end{align*}
with $\varepsilon=\textup{sign}(N)$.

Next we come to the Dirac operator.
With the $2\times 2$ Pauli matrices
$$
\sigma_+:=\begin{pmatrix}0&1\\0&0\end{pmatrix},\qquad
\sigma_0:=\begin{pmatrix}1&0\\0&-1\end{pmatrix},\qquad
\sigma_-:=\begin{pmatrix}0&0\\1&0\end{pmatrix},
$$ one has the relations
\begin{align}\label{pauli1}
\sigma_+\sigma_-=\begin{pmatrix}1&0\\0&0\end{pmatrix},\qquad\sigma_0^2=\begin{pmatrix}1&0\\0&1\end{pmatrix}, \qquad
\sigma_-\sigma_+=\begin{pmatrix}0&0\\0&1\end{pmatrix}, \\
\sigma_0\sigma_+=\sigma_+,\quad \sigma_+\sigma_0=-\sigma_+,\quad \sigma_+^2=\sigma_-^2=0,\quad
\sigma_-\sigma_0=\sigma_-,\quad
\sigma_0\sigma_-=-\sigma_- . \nn
\end{align}
Further, we use the differential operators $D_\pm,$ $D_0$,
\begin{equation}\label{diracbits}D_\pm:=L_\pm, \qquad
D_0:=L_0+q^{-2}L_z=\iq(q-\iq)^2(C_q+\tfrac{1}{4}-[\tfrac{1}{2}]^2),
\end{equation}
having used the expression \eqref{LCas} for the last equality.  As will be clearly momentarily, the use of
$D_0$ instead of $L_0$ (the extra $L_z$ vanishing identically on $\A[S^2_q]$)
will lead to a Dirac operator whose square is diagonal.
We define a
Dirac operator $D:\mathcal{S}\to \mathcal{S}$ by
\begin{equation}\label{dirac}D=D_+\sigma_++D_0\sigma_0+D_-\sigma_-,\end{equation}
where the $2\times 2$ Pauli matrices $\sigma_\pm$, $\sigma_0$ act
upon the column vector of $\mathcal{S}$ by left multiplication and
the vector fields $D_\pm$, $D_0$ operate {\em via} the left action
of $\usu$ (using the symbol $\tr$, which we omit from now on). As
mentioned above, elements $a\in\A[S^2_q]$ act as multiplicative
operators on $\mathcal{S}$ via the representations $\pi_{\pm}$:
$$\pi:\A[S^2_q]\to \textup{End}(\mathcal{S}),\qquad
\pi(a):=\begin{pmatrix}\pi_{+}(a)&0\\0&\pi_{-}(a)\end{pmatrix}$$
although we will not always explicitly denote the representation
$\pi$.

\begin{prop}
The Dirac operator $D:\mathcal{S}\to\mathcal{S}$ obeys
$$[D,a]=(L_+ a)\sigma_+ +(L_0 a)\sigma_0+(L_-
a)\sigma_-$$ for each $a\in \A[S^2_q]$.
\end{prop}
\proof
For $\psi=\begin{pmatrix}\psi_+&\psi_-\end{pmatrix}^{\textup{tr}}\in\mathcal{S}_+\oplus\mathcal{S}_-$,
using the derivation property of the vector fields $D_\pm$, $D_0$,
the commutator $[D,a]$ works out to be
\begin{align*}[D,a]\psi&=\begin{pmatrix}(D_+a)\psi_-\\0\end{pmatrix}+\begin{pmatrix}(D_0a)\psi_+\\-(D_0a)\psi_-\end{pmatrix}
+\begin{pmatrix}0\\(D_-a)\psi_+\end{pmatrix}\\ &=\left((D_+
a)\sigma_++(D_0 a)\sigma_0+(D_- a)\sigma_-\right)\psi.\end{align*}
To obtain the desired result, one simply substitutes $D_\pm=L_\pm$
and $D_0=L_0+q^{-2}L_z$, observing that $L_za=0$ for all $a\in
\A[S^2_q]$.
\endproof
\noindent This also shows that for all $a\in\A[S^2_q]$ the
commutator $[D,a]$ recovers the one-form $\D a$, acting on the
spinors $\mathcal{S}$ by `Clifford multiplication'.

The summand $D_+\sigma_+ +D_-\sigma_-$ in the operator \eqref{dirac}
is precisely the Dirac operator of \cite{ds:est}, corresponding
\cite{sw04} to the `two-dimensional' differential calculus on the
sphere $S^2_q$. The extra term $D_0$ in our Dirac operator is the
origin of the extra `direction' in the calculus $\Omega^1S^2_q$. It
is clear from~\eqref{LCas} that $D_0$ vanishes when $q\to 1$, whence
the classical limit of our construction is just the canonical
spectral triple on the classical two-sphere $S^2$.

Next, we compute the spectrum of the
Dirac operator. We shall use the identities
\begin{align}\label{qcomm1}
L_+L_-=qEFK^{-2}&=q\left(C_q+\tfrac{1}{4}-\frac{\iq K^2-2+qK^{-2}}{(q-\iq)^2}\right)K^{-2},\\
L_-L_+=\iq FEK^{-2}&=\iq\left(C_q+\tfrac{1}{4}-\frac{qK^{2}-2+\iq
K^{-2}}{(q-\iq)^2}\right)K^{-2}, \nn
\end{align} each obtained using the
expression \eqref{casimir} for the quantum Casimir $C_q$. Moreover,
we know from~\eqref{leftgens1} that for all $\psi_\pm\in
\mathcal{S}_{\pm}$ we have \begin{equation}\label{K
act1}K^2\psi_\pm=q^{\pm 1}\psi_\pm,\qquad K^{-2}\psi_\pm=q^{\mp
1}\psi_\pm.\end{equation}

These facts lead to the following result.
\begin{prop}\label{pr:sqd}
The Dirac operator $D$ obeys
$$
D^2=q^{-2}\nu^4 \left((C_q+\tfrac{1}{4}-[\tfrac{1}{2}]^2)\right)^2+\left(C_q+\tfrac{1}{4}\right),
$$
where $C_q$ is the quantum Casimir.
\end{prop}
\proof Using the Pauli relations \eqref{pauli1} one
computes that, for
$\psi=\begin{pmatrix}\psi_+&\psi_-\end{pmatrix}^{\textrm{tr}}\in
\mathcal{S}$,
\begin{align}\label{sqdirac}D^2\psi=&D_0^2\begin{pmatrix}1&0\\0&1\end{pmatrix}\psi+D_+D_-\begin{pmatrix}1&0\\0&0\end{pmatrix}\psi+D_-D_+\begin{pmatrix}0&0\\0&1\end{pmatrix}\psi .
\end{align}
The crucial fact in this calculation is that $D_0$ is a function of
the Casimir $C_q$ and therefore commutes with $D_\pm$. Next, using
the relations \eqref{qcomm1} and~\eqref{K act1} we find
\begin{equation*}\label{cas}D_{\pm}D_{\mp}\psi_{\pm}=\left(
C_q+\tfrac{1}{4}\right)\psi_{\pm}\end{equation*} for each
$\psi_\pm\in\mathcal{S}_\pm$.
Furthermore, we have that
$$D^2_0=q^{-2}\nu^4\left((C_q+\tfrac{1}{4}-[\tfrac{1}{2}]^2)\right)^2.$$
Substituting these expressions into~\eqref{sqdirac} yields the formula as claimed.
\endproof

As an immediate consequence we obtain the spectrum of our Dirac
operator $D$.
\begin{cor}\label{co:dspec}
The Dirac operator $D$ defined in \eqref{dirac} has spectrum
$$
\textup{Spec}(D)=\left\{ \pm\left( q^{-2}\nu^4 \, [j]^2[j+1]^2+[j+\tfrac{1}{2}]^2 \right)^{1/2}
~|~j\in\mathbb{N}+\tfrac{1}{2} \right\}
$$
with multiplicities $2j+1$.
\end{cor}

\proof The eigenvalues of $C_q$ are given in~\eqref{cas-spec}:
each $|j,m\ra_\pm$ is an eigenvector with eigenvalue
$[j+\tfrac{1}{2}]^2-\tfrac{1}{4}$, whence the multiplicity of the
$j$-th eigenvalue is $2(2j+1)$. From the expression for $D^2$ in
Prop.~\ref{pr:sqd}, we read off its eigenvalues using
those for $C_q$, yielding
\begin{equation}\textrm{Spec}(D^2)
=\left\{\lambda_j:= q^{-2}\nu^4 [j]^2 [j+1]^2 +[j+\tfrac{1}{2}]^2~|~j\in
\mathbb{N}+\tfrac{1}{2}\right\},
\end{equation} each having
multiplicity $2(2j+1)$. Here we have used the identity
$[j+\tfrac{1}{2}]^2-[\tfrac{1}{2}]^2=[j][j+1]$.
The eigenvalues of $D$ are therefore just
$\pm \lambda_j^{1/2} $ with multiplicities $2j+1$. \endproof

By inspection, we see that the eigenvalues of $|D|$ grow not
faster than $q^{-2j}$ for large $j$, in contrast with the Dirac
operator of \cite{ds:est}, whose eigenvalues diverge not faster than
$q^{-j}$. It is the extra term $D_0$ which
accounts for this behaviour.

This result immediately gives us an expression for $D$ in terms of
an orthonormal basis of eigenspinors $\vu$, $\vd$ defined by
\begin{equation}\label{diagdirac}
D|j,m;\uparrow\ra=\mu_j|j,m;\uparrow\ra,\quad D
|j,m;\downarrow\ra=-\mu_j|j,m;\downarrow\ra
\end{equation}
with eigenvalues
$$
\mu_j:=\left( q^{-2}\nu^4[j]^2[j+1]^2 +[j+\tfrac{1}{2}]^2\right)^{1/2}.
$$
To proceed further, it will be necessary to have an explicit
description of these eigenspinors in terms of the basic spinors
$|j,m\ra_\pm$. By evaluating the actions of $D_\pm$, $D_0$ on
$\mathcal{S}$ one finds that the Dirac operator is
\begin{equation}\label{basisdirac}
D|j,m\ra_\pm=\pm\iq\nu^2[j][j+1]|j,m\ra_\pm+[j+\tfrac{1}{2}]|j,m\ra_\mp,\end{equation}
the first term corresponding to the action of $D_0\sigma_0$, the
second to the action of $D_\pm\sigma_\pm$. Knowing the eigenvalues
of $D$, we find the corresponding eigenspinors to be
\begin{align}
\vu&:=\tfrac{1}{\sqrt{2\mu_j}}\left(-\zeta^+_j|j,m\ra_+-\zeta^-_j|j,m\ra_-\right),\\
\vd&:=\tfrac{1}{\sqrt{2\mu_j}}\left(-\zeta^-_j|j,m\ra_++\zeta^+_j|j,m\ra_-\right),
\nn
\end{align}
for $m=-j,-j+1,\ldots,j-1,j$ and $j\in\mathbb{N}+\tfrac{1}{2}$,
where we have written
\begin{equation}\label{zeta}
\zeta^+_j=\sqrt{\mu_j+\iq\nu^2[j][j+1]},\qquad
\zeta^-_j=\sqrt{\mu_j-\iq\nu^2[j][j+1]}.
\end{equation}
On the two-dimensional subspace $V_{j,m}$ spanned by $|j,m\ra_+,
|j,m\ra_-$ for fixed values of $j,m$, the operator which
diagonalises $D$ is just the orthogonal matrix
\begin{equation}\label{W}
W_{j}:=\frac{1}{\sqrt{2\mu_j}}\begin{pmatrix}-\zeta^+_j&-\zeta^-_j\\-\zeta^-_j&\zeta^+_j\end{pmatrix}.
\end{equation}
We write $W:\mathcal{H}\to\mathcal{H}$ for the closure of the
operator defined by the matrices $W_j$,
$j\in\mathbb{N}+\tfrac{1}{2}$.

\subsection{Spectral properties of $S^2_q$} We now show that the
datum $(\A[S^2_q],\mathcal{H},D)$ fulfils the conditions required of
a spectral triple, which we then equip with a real structure in the
sense of Definition~\ref{def:real}.

\begin{thm}\label{th:spec}The datum
$(\A(S^2_q),\mathcal{H},D)$ constitutes a unital spectral triple
over the sphere $S^2_q$ with metric dimension zero.\end{thm}

\proof For each $a\in \A[S^2_q]$ the commutator $[D,a]$ acts on
$\mathcal{S}$ by multiplication operators and is therefore itself a
bounded operator. In fact, for the summand $D_+\sigma_+
+D_-\sigma_-$ this goes as in \cite{ds:est}, whereas for the term
$D_0$ one gets multiplication by $L_0 a$ which belongs to
$\A[S^2_q]$ itself. The operator $D$ clearly satisfies $D=D^*$ on
the dense domain $\mathcal{S}$ of $\mathcal{H}$. From
Cor.~\ref{co:dspec} it is clear that the only accumulation points of
the spectrum of $D$ are at infinity, so the resolvent of $D$ is
compact. Since the eigenvalues of $D$ grow exponentially with
$j\in\mathbb{N}+\tfrac{1}{2}$, the metric dimension is just
zero.\endproof

\begin{prop}\label{pr:grading} With the $\ZZ_2$-grading
$\Gamma:\mathcal{H}\to \mathcal{H}$ defined by
$$
\Gamma \vu:=\vd,\qquad \Gamma\vd:=\vu
$$
on the orthonormal basis
\eqref{diagdirac} and extended by $\A[S^2_q]$-linearity, the datum
$(\A[S^2_q],\mathcal{H},D,\Gamma)$ constitutes an even spectral
triple.\end{prop}

\proof It is obvious that $\Gamma^2=1$ and $\Gamma=\Gamma^*$. The
property $\Gamma D+D\Gamma=0$ follows from the fact that $\Gamma$
interchanges the $+\mu_j$ and $-\mu_j$ eigenspaces of $D$, as may be
verified directly on the basis vectors \eqref{diagdirac}.\endproof

Next a real structure. Since we have made the same choice for the
spinors as in \cite{ds:est}, it is tempting to take the same real
structure as well. However, one quickly finds that this choice is
unsuitable, since it neither commutes nor anti-commutes with our
Dirac operator $D$. The reason for this lies mainly in the fact that
the term $D_0$ in our Dirac operator \eqref{dirac} is proportional
to the Casimir operator, which is rather a `second order
differential operator', if anything. Instead, we define an
anti-unitary operator $J:\mathcal{H}\to\mathcal{H}$ in terms of its
action on the orthonormal basis \eqref{diagdirac} by
$$
J\vu=(-1)^{m+1/2}|j,-m;\uparrow\ra,\qquad
J\vd=(-1)^{m+1/2}|j,-m;\downarrow\ra
$$
and seek to show that this $J$ equips the datum
$(\A[S^2_q],\mathcal{H},D,\Gamma)$ with a real structure.
It is not difficult to check that the $J$ above is equivariant under the right action
of $\usu$ on $\mathcal{H}$, making it a particularly natural choice.

\begin{prop}\label{pr:KO} The operator $J$ satisfies
$J^2=-1$, $DJ=JD$ and $\Gamma J=-J\Gamma$.\end{prop}

\proof The fact that $J^2=-1$ is immediate. We find that
\begin{align*}(DJ-JD)\vu&=(-1)^{m+1/2}D|j,-m;\uparrow\ra-\mu_jD\vu\\
&=(-1)^{m+1/2}\mu_j|j,-m;\uparrow\ra-(-1)^{m+1/2}\mu_j|j,-m;\uparrow\ra=0,\\
(J\Gamma+\Gamma J)\vu&=J\vd-(-1)^{m+1/2}\Gamma|j,-m;\uparrow\ra\\
&=(-1)^{m+1/2}|j,-m;\downarrow\ra-(-1)^{m+1/2}|j,-m;\downarrow\ra=0,
\end{align*}
where we have used anti-linearity of $J$.
Similar computations hold on $\vd$.
\endproof

Aiming at (modified) commutant and first order conditions as in
Definition~\ref{def:real}, and having in mind the strategy of
\cite{dalw}, we denote by $L_q$ the positive trace-class operator
defined by
$$
L_q|j,m\ra_\pm:=q^j|j,m\ra_\pm , \qquad
j\in\mathbb{N}+\tfrac{1}{2} ,
$$
on $\mathcal{H}$  and let $\mathcal{K}_q$ be the
two-sided ideal of $\mathcal{B}(\mathcal{H})$ generated by the
operators $L_q$. The ideal $\mathcal{K}_q$ is an ideal of
infinitesimals of arbitrarily high order and so we take
$\mathscr{I}=\mathcal{K}_q$ as our operator ideal in
Definition~\ref{def:real}. Thus, to prove that $J$ defines a real
structure, it remains to check that the commutant property and first
order condition in~\eqref{infsmod} are satisfied.

The strategy of \cite{dalw} is based on the fact that the operators $\pi(x_i)$, $i=-1,0,1$, can be
`approximated' by operators acting diagonally on the Hilbert space of spinors.
Specifically, these
operators $z_i$, $i=-1,0,1$, on $\mathcal{H}$ are defined by
\begin{multline}\label{approxops}
z_i|j,m\ra_\pm
=\alpha^-_i(j,m;0)|j-1,m+i\ra_\pm+\alpha^0_i(j,m;0)|j,m+i\ra_\pm+\alpha^+_i(j,m;0)|j+1,m+i\ra_\pm.
\end{multline}
The coefficients are exactly the ones used in~\eqref{irreps},
unless $|m+i|>j+\nu$ for $\nu=-1,0,1$, in which case we set
$\alpha^\nu_i(j,m;0)=0$. Momentarily we shall show that the
operators $z_i$ approximate the operators $\pi(x_i)$ modulo the
ideal $\mathcal{K}_q$, but to do this we first need the following
technical lemma.

\begin{lem}\label{le:bound}
With $W_j$, $j\in\mathbb{N}+\tfrac{1}{2}$, the operators in~\eqref{W}, there exists a constant $C$ (independent of $j$) such
that
$$
||W_jW^*_{j+1}-1||<Cq^{j}
$$
for all $j\in \mathbb{N}+\tfrac{1}{2}$.\end{lem}

\proof One evaluates the norm $||W_jW^*_{j+1}-1||$ by computing the
eigenvalues of the $2\times 2$ matrix $W_jW^*_{j+1}-1$ and choosing
the larger of the two, finding it to be
$$
||W_jW^*_{j+1}-1||=\frac{\zeta^+_{j}\zeta^+_{j+1}+\zeta^-_{j}\zeta^-_{j+1}-\zeta^-_{j}\zeta^+_{j+1}+\zeta^+_{j}\zeta^-_{j+1}-2\sqrt{\mu_j\mu_{j+1}}}{2\sqrt{\mu_j\mu_{j+1}}}.
$$
Using the inequalities $[j]<(q-\iq)^{-1}q^{-j}$ and
$[j]^{-1}<q^{j-1}$, elementary estimates for each of the terms in
this expression yield that $\zeta^\pm_j<C'q^{-j}$ and
$\sqrt{\mu_j\mu_{j+1}}<C''q^{-2j}$ for real constants $C',C''$, so
it appears at first glance that the above norm has an $O(1)$
behaviour. However, a more detailed analysis shows that the
coefficient of $q^{-2j}$ in the numerator is in fact zero; the
behaviour of the numerator is therefore $O(q^{-j})$ and we have our
result.\endproof


\begin{prop}\label{pr:approx} There exist bounded
operators $A_i$, $B_i$, $i=-1,0,1$, such that
$$
\pi(x_i)-z_i=A_iL_q=L_qB_i
$$
when acting upon the basis vectors
$|j,m;\uparrow\downarrow\ra$. In particular,
$\pi(x_i)-z_i\in\mathcal{K}_q$ for $i=-1,0,1$.
\end{prop}

\proof From \cite[Lem.~4.4]{dalw}, there exist bounded
operators $A_i$, $B_i$, $i=-1,0,1$ such that
$$
\pi(x_i)-z_i=A_iL_q=L_qB_i
$$
with respect to the basis $|j,m\ra_\pm$ of $\mathcal{H}$, and so the
operators $\pi(x_i)$ are approximated by the operators $z_i$ modulo
the ideal $\mathcal{K}_q$ of infinitesimals. We need to check that
using the operator $W$ to change the basis vectors from
$|j,m\ra_\pm$ to $|j,m;\uparrow\downarrow\ra$ does not spoil this
approximation property. Evaluating $W_jz_iW_j^*-z_i$ on
$|j,m;\uparrow\downarrow\ra$ gives
\begin{multline*}
(W_jz_iW_j^*-z_i)|j,m;\uparrow\downarrow\ra =\alpha^-_i(j,m;0)(W_{j-1}W_j^*-1)|j-1,m+i;\uparrow\downarrow\ra\\
 +\alpha^+_i(j,m;0)(W_jW_{j+1}^*-1)|j+1,m+i;\uparrow\downarrow\ra.
\end{multline*}
This and Lemma~\ref{le:bound} yield that
$W_jz_iW^*_j-z_i\in\mathcal{K}_q$ for all $i=-1,0,1$ and all
$j\in\mathbb{N}+\tfrac{1}{2}$.\endproof

As a consequence, we immediately get the commutant property, the
first of the two conditions in~\eqref{infsmod}.

\begin{prop}\label{pr:comm} For all $a,b\in\A[S^2_q]$ we have
$[\pi(a),J\pi(b)J^{-1}] \in\mathcal{K}_q$.\end{prop}

\proof From the derivation property of commutators, it suffices to
check this only for the generators $x_{-1},x_0,x_1$ of $\A[S^2_q]$.
With the operators $z_{-1},z_0,z_1$ defined in~\eqref{approxops}, we have
\begin{multline}\label{z-action}
Jz_kJ^{-1}|j,m\ra_\pm  =(-1)^k\left(\alpha^-_k(j,-m;0)|j-1,m-k\ra_\pm\right. \\
 \left.+\alpha^0_k(j,-m;0)|j,m-k\ra_\pm+\alpha^+_k(j,-m;0)|j+1,m-k\ra_\pm\right).
\end{multline}
Using this, one computes as in \cite[Lem.~6.2]{dalw} that
\begin{equation}\label{commut}
[z_i,Jz_kJ^{-1}]=0,\qquad i,k=-1,0,1.
\end{equation}
It is straightforward to check that
$$
[\pi(x_i),J\pi(x_k)J^{-1}]=[\pi(x_i)-z_i,J\pi(x_k)J^{-1}]+[z_i,J\left(\pi(x_k)-z_k\right)J^{-1}]+[z_i,Jz_kJ^{-1}],
$$
whence the assertion follows from
Prop.~\ref{pr:approx}.\endproof

We are now ready for our main theorem regarding the differential
structure of $S^2_q$.

\begin{thm}The datum
$(\A(S^2_q),\mathcal{H},D,\Gamma,J)$ constitutes a real even unital
spectral triple (up to infinitesimals) with \textup{KO}-dimension
equal to two.\end{thm}

\proof Having already established Props.~\ref{pr:KO} and
\ref{pr:comm}, it remains to verify the first order condition for
$D$, namely that $[[D,a],JaJ^{-1}]\in\mathcal{K}_q$ for all
$a\in\A[S^2_q]$. For this, we split the Dirac operator into two
pieces, $D=D_\Delta+D_\Omega$, where $D_\Delta=D_0\sigma_0$ and
$D_\Omega=D_-\sigma_-+D_+\sigma_+$. By linearity
it suffices to check the first order condition for $D_\Delta$ and
$D_\Omega$ individually.

Since $D_0$ is a function of the Casimir, each $a\in\A[S^2_q]$ is an
eigenfunction for the derivation $[D_\Delta,\,\cdot\,]$, whence the
first order condition for $D_\Delta$ follows immediately from the
commutant property in Prop.~\ref{pr:comm}. On the other hand, the
component $D_\Omega$ has eigenvalues $\pm \gamma_j$,
$\gamma_j:=[j+\tfrac{1}{2}]$, whose growth with $j$ obeys
$\gamma_j<Cq^{-j}$ for $C$ a real constant (as already mentioned,
$D_\Omega$ is precisely the Dirac operator considered in
\cite{ds:est}). It is easy to compute that
\begin{multline*}
[D_\Omega,z_i]|j,m\ra_\pm =(\gamma_{j-1}-\gamma_j)\alpha^-_i(j,m;0)|j-1,m+i\ra_\mp\\
 \qquad\qquad+(\gamma_{j+1}-\gamma_j)\alpha^+_i(j,m;0)|j+1,m+i\ra_\mp.
\end{multline*}
Using this expression, together with~\eqref{z-action}, one calculates the action of the commutators
$[[D_\Omega,z_i],Jz_kJ^{-1}]$ for $i,k=-1,0,1$ and finds them to be
a sum of five independent weighted shift operators with weights
$S^\nu_{i,k}(j,m)$, $\nu=-2,\ldots,2$, {\em i.e.}
$$
[[D_\Omega,z_i],Jz_kJ^{-1}]|j,m\ra_\pm=\sum_{\nu=-2}^2
S^\nu_{i,k}(j,m)|j+\nu,m+i-k\ra_\pm.
$$
These weights $S^\nu_{i,k}(j,m)$ are estimated using exactly the
same method as in \cite[Prop.~6.5]{dalw}. In our case, the growth
condition for $\gamma_j$ is sufficient to guarantee that
$|S^\nu_{i,k}(j,m)|<C'q^j$ for some real constant $C'$. We conclude
that $[[D_\Omega,z_i],Jz_kJ^{-1}]\in\mathcal{K}_q$ for all
$i,k=-1,0,1$. Since the $z_i$ approximate the operators $\pi(x_i)$
modulo $\mathcal{K}_q$, the proof is complete.\endproof

\subsection*{Acknowledgments}
Both authors were partially supported by the Italian Project
`Cofin08--Noncommutative Geometry, Quantum Groups and Applications'.
SB is grateful to INdAM--GNSAGA for support and the Department of
Mathematics at the University of Trieste for its hospitality. We
thank Francesco D'Andrea for very useful comments.

\end{document}